\documentclass{gtpart}
\usepackage{amsfonts}
\usepackage{amssymb}
\usepackage{amsthm}
\usepackage{verbatim}
\usepackage[dvips]{graphicx}

\author{Katharine C. Walker}
\givenname{Katharine}
\surname{Walker}
\address{Department of Mathematics\\University of Michigan\\3839 East Hall\\Ann Arbor, MI 48103}

\email{kaceyw@umich.edu}

\theoremstyle{plain}
\newtheorem{theorem}{Theorem}[section]
\newtheorem{lemma}[theorem]{Lemma}
\newtheorem{definition}[theorem]{Definition}
\newtheorem{corollary}[theorem]{Corollary}
\newtheorem{prop}[theorem]{Proposition}

\newtheorem{constr}[theorem]{Construction}

\renewcommand{\C}{\mathbb{C}}
\renewcommand{\Z}{\mathbb{Z}}
\renewcommand{\R}{\mathbb{R}}
\newcommand{\N}{\mathbb{N}}
\renewcommand{\P}{\mathbb{P}}
\newcommand{\bs}{\backslash}
\renewcommand{\epsilon}{\varepsilon}

\renewcommand{\phi}{\varphi}

\newcommand{\into}{\hookrightarrow}
\renewcommand{\Q}{\mathcal{Q}}
\newcommand{\bQ}{\overline{\mathcal{Q}}}

\newcommand{\M}{\mathcal{M}}

\newcommand{\tM}{\tilde{M}}
\newcommand{\tq}{\tilde{q}}
\newcommand{\tp}{\tilde{p}}
\newcommand{\hM}{\hat{M}}
\newcommand{\hq}{\hat{q}}
\newcommand{\he}{\hat{e}}
\newcommand{\hh}{\hat{h}}
\newcommand{\T}{\mathcal{T}}

\renewcommand{\l}{\lambda}
\renewcommand{\a}{\alpha}

\newcommand{\g}{\gamma}

\newcommand{\ep}{\epsilon}

\begin{document}
\title[Quotient groups of the fundamental groups]{Quotient groups of the fundamental groups of certain strata of the moduli space of quadratic differentials.} 
\begin{abstract} In this paper, we study fundamental groups of strata of the moduli space of quadratic differentials. We use certain properties of the Abel-Jacobi map, combined with local surgeries on quadratic differentials, to construct quotient groups of the fundamental groups for a particular family of strata.  
\end{abstract}
\maketitle
\section{Introduction}
In \cite{KZ1}, Kontsevich and Zorich conjecture that the fundamental groups of strata of the moduli space of abelian differentials are commensurable with various mapping class groups. In this paper we consider a similar question for strata of quadratic differentials, and in particular we construct a quotient group of the fundamental group for a certain family of strata. We do so by mapping a stratum into a larger configuration space of points on surfaces,  and showing that the image of the fundamental group of the stratum under this map is in the kernel of a version of the Abel-Jacobi map. We then construct a set of generators for the kernel of the Abel-Jacobi map, and show that in some cases the image of the fundamental group of the stratum in the fundamental group of the configuration space is equal to this kernel. 

More specifically, let $\Q_g$ be the space of quadratic differentials over Teichmuller space, $\T_g$, and let  $\l=(k_1,...,k_n)$ be a partition of $4g-4$. Define $\Q_g(k_1,...,k_n)=\Q_\l$ to be the subset of $\Q_g$ of quadratic differentials with $n$ zeros of order $k_1,...,k_n$. Let $\bQ_g$ and $\bQ_\l$ be the analogous spaces over moduli space, $\M_g$. We are interested in $\pi_1(\bQ_\l)$; however, when $\Q_\l$ and $\bQ_\l$ are both connected we have a short exact sequence $\pi_1(\Q_\l) \to \pi_1(\bQ_\l) \to \Gamma_g$, where $\Gamma_g$ is the genus $g$ mapping class group. In many cases (although not all) both $\Q_\l$ and $\bQ_\l$ are connected, so we focus on proving results about $\pi_1(\Q_\l)$. 

To do this we first embed $\Q_\l$ into a larger configuration space. In particular, to any partition, $\l$, we associate a generalized symmetric group, $S_\l$, that allows points of equal weights to be exchanged. For $M \in \T_g$ and $\l$ of length $n$, let $M^{[n]}$ denote the space of $n$ ordered distinct marked points on $M$, and let $Sym^\l(M)$ denote $M^{[n]} / S_\l$. Define $Sym_g^\l$ to be the associated bundle over $\T_g$, and $Pic_g^{4g-4}$ the bundle over $\T_g$ with fiber  $Pic^{4g-4}(M)$, the Picard variety parametrizing line bundles on $M$ of degree $4g-4$. Then we have the following maps:
\begin{equation} \label{eq:ses}
\Q_\l \stackrel{i}{\to} Sym_g^\l \stackrel{AJ}{\to} Pic_g^{4g-4}
\end{equation}
The first map is given by considering the zeroes of a quadratic differential as weighted marked points. The second map is the Abel-Jacobi map, given by mapping a divisor to its associated line bundle. 
 The maps in (\ref{eq:ses}) induce a sequence of maps: 
\begin{equation} \label{eq:pi1}
\pi_1(\Q_\l) \stackrel{i_*}{\to} \pi_1(Sym_g^\l) \stackrel{AJ_*}{\to} \pi_1(Pic_g^{4g-4})=H_1(M, \Z)
\end{equation}
We show that $AJ_* \circ i_*:\pi_1(\Q_\l) \to H_1(\Sigma, \Z)$ is trivial, so the image of $i_*$ will be in the kernel of $AJ_*$. In the case where there are at least $O(\sqrt{g})$ zeroes of order 1 in $\l$, we are able to construct a set of generators for the kernel of $AJ_*$.
\begin{theorem} \label{thm:pi1intro}
Let $\l=(k_1^a,k_2^{b_2},...,k_m^{b_m})$, with $\sum_{i=2}^m b_i=b$ and $a \ge \frac{3 + \sqrt{9 +8(2g+b-2)}}{2}$. Then the kernel of $AJ_*:\pi_1(Sym_g^\l) \to H_1(M, \Z)$ is generated by transpositions of zeroes of equal weight, squares of transpositions of zeroes of unequal weight, moving sets of points of equal weight opposite ways around generators of $\pi_1(M)$, and in some cases moving single points around homologically trivial curves. 
\end{theorem}

A more precise statement of the theorem is given in Section \ref{sec:ker}. To show that the elements detailed in Theorem \ref{thm:pi1intro} are contained in $\pi_1(\Q_\l)$ we first follow a variety of authors, including \cite{EMZ}, \cite{L1}, and \cite{KZ}, to create local surgeries on surfaces with quadratic differentials that affect the surface and quadratic differential in only a small area around a zero. We then take advantage of the fact that for any genus 0 stratum, $Q_\l \cong Sym^\l(\P^1)$ (in other words, a quadratic differential on $\P^1$ may have zeroes at any set of points) to create explicit curves of quadratic differentials in the hyperelliptic loci of certain strata. This leads to the following theorem: 
\begin{theorem} \label{thm:genintro}
 Let $\l=(1^a,k_1,...,k_n)$ with $a >$ max$\{g+5,k_1,..., k_n \}$, all $k_i$ even, and some $k_i=k_j$. Then $im(i_*:\pi_1(\Q_\l, (M,q)) \to \pi_1(Sym_g^\l))=ker(AJ_*)$.
\end{theorem}

Theorem \ref{thm:genintro} states that for certain $\l$ we can describe the image of $\pi_1(\Q_\l^0)$ in $\pi_1(Sym_g^\l)$. However, the kernel of $\pi_1(\Q_\l^0) \to \pi_1(Sym_g^\l)$ may be non-trivial, so we have created a quotient group of $\pi_1(\Q_\l^0)$.

The structure of the paper is as follows. In Section \ref{sec:pre} we give some general background. In Section \ref{sec:pi1} we collect some results about $\pi_1(Sym_g^\l)$, and in Section \ref{sec:ker} we construct a set of generators for the kernel of $AJ_*$ for $\l$ with sufficiently many zeroes of the same order. In Section \ref{sec:adj} we develop some local surgeries that allow us to construct elements in this kernel. In Sections \ref{sec:hyp} and \ref{sec:null} we use the results of Section \ref{sec:adj} to construct explicit elements in the $ker(AJ_*)$. Section \ref{sec:sum} summarizes when we have the image of $\pi_1(\Q_\l) \to \pi_1(Sym_g^\l)$ equal to the kernel of $AJ_*:\pi_1(Sym^\l_g) \to H_1(M, \Z)$, as well as some of the difficulties in analyzing the kernel of $i_*$. Theorems \ref{thm:pi1intro} and \ref{thm:genintro} are proved in Sections \ref{sec:ker} and \ref{sec:sum}, respectively. 

\section{Preliminary Definitions} \label{sec:pre}
A \textit{meromorphic quadratic differential}, $q$, on a Riemann surface, $M$, is a meromorphic section of the square of the canonical bundle, $K$, of $M$. In local coordinates $q$ assigns to each $(U_\alpha, z_\alpha)$ a meromorphic function $f_\alpha$ such that:
$$f_\beta(z_\beta)(\frac{dz_\beta}{dz_\alpha})^2 = f_\alpha(z_\alpha), dz_\beta = \frac{dz_\beta}{dz_\alpha} dz_\alpha$$ on $U_\alpha \cap U_\beta$. 

A \textit{horizontal trajectory}, or simply a trajectory, of a quadratic
differential, $q$, on $M$ is a smooth curve $\gamma:[0,1] \to M$ such that
$f(\gamma(t)) (\gamma'(t)dt)^2$ is real and positive for all $t$. Similarly, a vertical trajectory of $q$ is $\gamma:[0,1] \to M$ such that $(f(\gamma(t)) (\gamma'(t)dt)^2)$ is real negative, and a $\theta$-trajectory is $\gamma$ such that the argument of ($f(\gamma(t)) (\gamma'(t)dt)^2$ is $2\theta$.  (For $q=dz^2$ on $\C$ these correspond to straight lines of angle $\theta$.) Through every regular point of $q$ there exist unique horizontal and vertical trajectories, which are transverse. Near a zero of $q$ of order $n$,  $q = a_n z^n + a_{n+1}z^{n+1}+....$ and in fact by changing variables $\zeta = z(a_n+a_{n+1}z+...)^{1/n}$ we can define $q=\zeta^n$. Then the curves $\gamma(t)=te^{i\frac{2\pi}{(n+2)}k}$ for $k=0,1,....,n+1,$ are all horizontal trajectories and $n+2$ trajectories dead-end into a zero of order $n$. A trajectory between two critical points of $q$ is called a \textit{saddle connection}.

Any quadratic differential $q$ gives us a metric on $M$: 
$$|\gamma|_{q} = \int_{\gamma} |f(z)|^{1/2} |dz|$$ 
for $\gamma$ a real curve on $M$. Geodesics in this metric are unions of $\theta$-trajectories, with vertices at critcal points of $q$. In general the distance between two points in the $q$ metric is not well-defined because there will be a geodesic associated to each homology class of curves between the two points. However, given a choice of a disc or polygon containing two points the distance becomes well-defined. 

Through much of this paper we will be concerned not just with individual quadratic differentials but also their moduli spaces. Thus we consider the bundle, $\mathcal{E}$ over $\T_g$ with fiber $H^0(M, K^2)$, and define $\Q_g$ to be the subspace of the total space of $\mathcal{E}$ consisting of quadratic differentials that are not squares of Abelian differentials (sections that are not the squares of sections of $K$). Define $\Q_g(k_1,...,k_n)$ to be the subspace of $\Q_g$ of quadratic differentials with zeroes of order $k_1,...,k_n$, $\sum_1^n k_i=4g-4$ and $k_i \in \N$. The partitions of $4g-4$ give a natural stratification of $\Q_g$ and a single $\Q_g(k_1,...,k_n)$ is often called \textit{stratum} of $\Q_g$. 
We also occasionally consider an analog of $\mathcal{E}$ where each fiber is the space of meromorphic sections of $K^2$ with up to some fixed number of single poles, and then we can consider $\Q_g(k_1,...,k_n)$ where some of the $k_i=-1$. All of these spaces are well known to be manifolds (see \cite{V}, for example).  $\Q_g$ is also well know to be the cotangent bundle of $\T_g$, and is thus a complex manifold of dimension $6g-6$. Unless $n=1$, $\Q_g(k_1,...,k_n)$ is not closed in $\Q_g$ as zeroes may collide to form higher order zeroes.  
The mapping class group, $\Gamma_g$, acts on $\Q_g(k_1,...,k_n)$ by a
lift of its action on $\T_g$. 
We define the quotient of $\Q_g$ by $\Gamma_g$ to be $\bQ_g$, a space over $\M_g$, and similarly define $\bQ_g(k_1,....,k_n)$ as the quotient of $\Q_g(k_1,....,k_n)$ by $\Gamma_g$. Since $\Gamma_g$ does not act freely, the $\bQ_g(k_1,...,k_n)$ will be complex orbifolds.

In general we do not require the $k_i$ to be distinct, but if a stratum has multiple zeroes of the same order we will sometimes use the notation $\Q_g(k_1^n,k_2,...,k_n)$ to indicate $\Q_g(k_1,k_1,...,k_1,k_2,...,k_n)$. In general we will denote elements of $\Q_g(k_1,...,k_n)$ as $(M,q)$ where $M$ may be though of as a Riemann surface with the extra data of a homology basis attached. When the specific orders of the zeroes are not important we will sometimes let $\l=(k_1,...,k_n)$ denote a partition of $4g-4$ and let $\Q_\l=\Q_g(k_1,...,k_n)$. The \textit{length} of $\l$ will be the number of $k_i$ in the partition. The following sums up the structure of various strata: 
\begin{theorem} [Masur, Smillie, Veech]
Every $\bQ_g(k_1,...,k_n)$ is non-empty, with four exceptions:
$\bQ_1(\emptyset), \bQ_1(-1,1), \bQ_2(3,1), \bQ_2(4)$. With the exception of these four strata, the $\bQ_g(k_1,...,k_n)$ are complex orbifolds of dimension $2g-2+n$
\end{theorem}
The same is true of the $\Q_g(k_1,...,k_n)$, except that they are manifolds instead of orbifolds. Also:  

\begin{prop}\label{prp:sphere}
Any $\Q_0(k_1,...,k_n)$ is connected. 
\end{prop} 
Notice that if one has a (ramified) cover of some $M$, one can pull back back a quadratic differential on $M$ to get one on its cover. 

\begin{lemma} [Lanneau] \label{lem:ord} Let $\pi: \tilde{M} \rightarrow M$ be a ramified double
  cover, $q$ a quadratic differential on $M$ and $\tilde{q}$
  its pullback under $\pi$. Let $\tilde{p}$ be a ramification point of
  $\pi$ and $p=\pi(\tilde{p})$. Then, if $p$ is a singularity of order k of
  $q$, $\tilde{p}$ will be a singularity of order $2k+2$ on $\tilde{q}$.
  \end{lemma}  

One can see this by noticing that a singularity of order $k$ corresponds to a cone angle of $(k+2) \pi$; combining two points with cone angle $(k+2) \pi$ gives a cone angle of $2k+4 \pi$ or a singularity of order $2k+2$. 

Using double covers one may construct a continuous map between two different strata of quadratics differentials: 

\begin{constr}\label{con:lan} 
Let $\sum_1^n k_i = -4$, $n \ge 2g+2$. We can
  construct a local map $\Q_0(k_1,...,k_n) \to \Q_g(2k_1+2,...2k_{2g+2}+2,
  k_{2g+1}^2,...,k_n^2)$ by taking 2 copies of $(M',q') \in
  \Q_0(k_1,...,k_n)$, making $g+1$ cuts between the first $2g+2$
  marked points,and gluing along each of those cuts. This gives a
  surface $(M, q)$ of genus g such that each of the first $2g+2$
  zeroes of order $k_i$ goes to one with order $2k_i+2$, and we get 2 copies of each remaining zero. 
\end{constr}

Note that it is possible that the cover in Construction \ref{con:lan} will be the square of an Abelian differential, instead of a quadratic differential. However, if any singularity in the double cover is of odd order then it must be a quadratic differential.  

\begin{definition}
Define a quadratic differential $(M,q)\in \Q_g(k_1,...k_n)$ to be
\textit{hyperelliptic} if it is a double cover of some $(M',q') \in
\Q_0(k_1',...,k_m')$ as in Construction \ref{con:lan}. Define $\Q_g(k_1,...,k_n)$ to be \textit{hyperelliptic} if it contains
hyperelliptic quadratic differentials. 
\end{definition}

Finally, since $\bQ_\l \to \Q_\l \to \Gamma_g$ is a fibration, if $\Q_\l$ and $\bQ_\l$ are both connected then $\pi_1(\bQ_\l) \to \pi_1(\Q_\l) \to 
\pi_0(\Gamma_g)=\Gamma_g$ will be short exact. Interestingly, not all of the $\bQ_\l$ are actually connected. For those that are disconnected, the connected components are classified by whether or not hyperelliptic quadratic differentials form a full-dimensional subset of the stratum. Lanneau proves this in \cite{L1}.

\begin{theorem} [Lanneau]  \label{thm:lanneau}
For $g \ge 3$ the following strata have two connected components: 
\begin{enumerate}
\item $\bQ_g(4(g-k)-6, 4k+2)$, $k \ge 0, g-k \ge 2$
\item $\bQ_g((2(g-k)-3)^2, 4k+2)$, $k \ge 0, g-k \ge 1$
\item $\bQ_g((2(g-k)-3)^2, 2k+1^2)$, $k \ge 0, g-k \ge 2$
\end{enumerate}
and the rest have one component. 
For $g=0,1$ all strata are connected, and for $g=2$ $\bQ_2(3,3,-1,-1)$ and $\bQ_2(6,-1,-1)$ have two components, but all others are connected.
\end{theorem}

In \cite{W} we showed the following: 
\begin{theorem} \label{thm:c1} Let $m \ge g$.  Then any stratum of the form $\Q_g(1^{m},k_1^{n_1},...,k_l^{n_l})$ is connected.
\end{theorem} 

Then for $\l$ as in Theorem \ref{thm:c1} both $\Q_\l$ and $\bQ_\l$ are connected, and $\pi_1(\bQ_\l) \to \pi_1(\Q_\l) \to \Gamma_g$ is short exact.

\section{Surface Braid Groups} \label{sec:pi1}
In this section we collect some results about surface braid groups, to use in analyzing the kernel of $AJ_*$. 

Let $S_n$ be the standard symmetric group on $n$ letters. To any partition, $\l$, of $4g-4$ we associate a symmetric group, $S_\l$, which allows equal values to be exchanged. For example, to $(1^{4}, 2, 5^{2})$ we associate $S_4 \times S_2$. The \textit{length} of $\l$ will be the number of elements it contains. For a particular $M \in \T_g$ and partition $\l$ of length $n$, let $M^{[n]}$ denote the space of $n$ ordered distinct marked points on $M$, and let $Sym^\l(M)$ be $M^{[n]} / S_\l$. 

$\pi_1(M^{[n]})$ is well-known as the pure or special braid group of $n$ elements on a genus $g$ surface, which we will denote by $SB_n(M)$ or simply $SB_n$. Similarly, $\pi_1(M^{[n]}/S_n)$ is the full braid group on $M$, $B_n(M)$ or $B_n$. The generators of both $SB_n$ and $B_n$ are well-known. In particular, let $l_1,...,l_{2g}$ be $2g$ standard generators of $\pi_1(M)$ and let $\rho_{ij}$, $1 \le i \le n$, $1 \le j \le 2g$, denote an element of $SB_n$ such that $p_i$ follows a path that is homotopic to $l_j$. Let $D$ be any disk containing $p_k, p_l$, $1 \le k<l \le n$, and let $\kappa_{kl}$ be either generator of $\pi_1(D^{[2]}, (p_k,p_l))$. (We may extend this to an element of $\pi_1(M^{[n]})$ by letting the other $n-2$ points move along constant paths.)  Similarly, let $\sigma_{st}$ be either generator of $\pi_1(D^{[2]}/S_2, (p_s,p_t))$. The following theorem is classical. 

\begin{theorem} \label{thm:sct}
 $SB_n$ is generated by the $\rho_{ij}$, $1 \le i \le n$, $1 \le j \le 2g$, and the $\kappa_{kl}$, $1 \le k < l \le n$. $B_n$ is generated by the $\rho_{ij}$ and the $\sigma_{s(s+1)}$, $1 \le s < n$. 
\end{theorem}

It should be noted that there are multiple non-equivalent ways to define each of the $\rho_{ij}$, $\kappa_{kl}$, and $\sigma_{s(s+1)}$; however, any choice yields a generating set. 

Similarly, let $M_m$ denote $M - \{ p_{n+1},...,p_{n+m}\}$, where the $p_{n+i}$ are any distinct points on $M$, and let $SB_{n,m}$ denote $SB_n(M_m)$ and $B_{n,m}$ denote $B_n(M_m)$. (We will only be concerned with the topology of $M_m$, which does not depend on the choice of $\{p_{n+1},...,p_{n+m} \}$.) Let $\kappa_{kl}$, $1 \le k \le n$, $n < l \le n+m$, denote $p_k$ moving in a simple loop around $p_{l}$. Again the generators of both $SB_{n,m}$ and $B_{n,m}$ are well-known.

\begin{theorem} \label{thm:braid}
$SB_{n,m}$ is generated by the $\rho_{ij}$ and the $\kappa_{kl}$, $1 \le j \le 2g$, $1 \le i, k \le n$, $1 \le l \le n+m$. $B_{n,m}$ is generated by the $\rho_{ij}$, the $\sigma_{s(s+1)}$, $1 \le s < n$, and the $\kappa_{kl}$, $1 \le k \le n$, $n < l \le n+m$.  
\end{theorem}

We will primarily be interested in $\pi_1(Sym^\l(M))$, which we will denote by $B_\l(M)$ or $B_\l$. For $\l=(k^n)$, $B_\l$ is just $B_n$. For more complicated $\l$ we note that the covering map $M^{[n]} \to Sym^\l(M)$ is normal and thus $$SB_n \to B_\l \to S_\l$$ is a short exact sequence. This tells us the generators of $B_\l$:  

\begin{prop} \label{thm:prod}
Let $\l$ be a partition of $4g-4$ of length $n$. Then
$B_\l$ is generated by the $\rho_{ij}$, $1 \le i \le n$,
$1 \le j \le 2g$, and for each pair $1 \le k < l \le n$, either $\sigma_{kl}$ if $p_k$ and $p_l$ are of the same weight, or $\kappa_{kl}$ if $p_k$ and $p_l$ are of different weights. 
\end{prop}

In fact it is possible to generate $B_\l$ with fewer transpositions; however, this generating set will suffice for our purposes. Although it will not be explicitly used in this paper, it is also worth noting that for particular choices of generating sets, the relations among the generators of $SB_m, SB_{n,m}$, and $B_n$ are well-known, and thus the same will be true for any of their subgroups. 

Another classical result about surface braid groups is the following. Let $M^{[n]+[r]}$ be the space of all n-tuples, r-tuples of distinct ordered points that are disjoint. Then we have the following theorem from \cite{FN}:
 
\begin{theorem} [Fadell, Neuwirth] \label{thm:fn}
$(M_m)^{[n-r]+[r]} \to (M_m)^{[r]}$ is a fibration, with fiber $(M_{m+r})^{[n-r]}$.
\end{theorem}

This fibration induces a long exact sequence of homotopy groups. All higher homotopy groups are trivial, so
$$SB_{n-r,m+r}(M) \to SB_{n,m}(M) \to SB_{r,m}(M)$$ 
is short exact. 

Finally, we will sometimes want to distinguish the many different transpositions of two points on a surface. Let $p,p' \in M$ and and define an \textit{edge}, $e$, to be an embedding of the interval $[0,1]$ in $M$ with endpoints $p$ and $p'$.  Let $U \subset M$ be a contractible neighborhood of $e$. Then we define $\sigma_e$ to be either of the two generators of $\pi_1(U^{[2]}/S_2, (p,p'))$, with $\sigma_e^{-1}$ its inverse, and $\kappa_e$ to be either of the two generators of $\pi_1(U^{[2]}, (p,p'))$. For a particular $(M,q)$, let $P=\{p_1,....,p_n\}$ be the zeroes of $q$ and define $\bar{E}_{M,q}=\{e:I \into M | e(I) \cap P= e(0) \cup e(1)\}$ to be the set of all edges on $(M,q)$. Put the following equivalence relation on edges: $e \sim e'$ if there exists $h:I \times I \into M$ such that $h(I \times I) \cap P = \{0\} \times I \cup \{1\} \times I, h( \cdot, 0) = e, h(\cdot, 1)=e'\}$. We say $e \sim e'$ if the above is satisfied, and define $E_{M,q}=\bar{E}_{M,q}/\sim$. We will index the set of all transpositions associated to $(M,q)$ by $E_{M,q}$.

\section{The kernel of the Abel-Jacobi map} \label{sec:ker}
In this section we first define the Abel-Jacobi map, and then construct a set of generators for its kernel. 

Let $Sym^\l(M)$ denote $M^{[n]} / S_\l$, and define $Sym_g^\l$ to be bundle over $\T_g$ with fiber $Sym^\l(M)$. Similarly let $Pic^{4g-4}(M)$ be the Picard variety parametrizing line bunldes of degree $4g-4$, and define $Pic_g^{4g-4}$ to be the bundle over $\T_g$ with fiber  $Pic^{4g-4}(M)$. Let $\Lambda_g$ be the set of all partitions of $4g-4$ and define $Sym_g := \bigcup_{\l \in \Lambda_g} Sym_g^\l$.

 Then we have a sequence of maps: 
\begin{equation} \label{eq:seq}
\Q_g \stackrel{i}{\to} Sym_g \stackrel{AJ}{\to} Pic_g^{4g-4}.
\end{equation}
The first map is given by considering the
zeroes of a quadratic differential as marked points. The second map is the Abel-Jacobi map, given by mapping a divisor to its associated line bundle. The
composition of these maps is fiber-wise trivial because every element of
$\Q_g$ over a particular $M \in \T_g$ maps to $K_M^2 \in Pic_g^{4g-4}$. Since $\T_g$ is simply connected, the image of $\Q_g$ under the composition of these two maps is also simply connected. 

The maps in (\ref{eq:seq}) induce a sequence of maps: 
\begin{equation} 
\pi_1(\Q_g) \stackrel{i_*}{\to} \pi_1(Sym_g) \stackrel{AJ_*}{\to} \pi_1(Pic_g^{4g-4}) \cong H_1(M, \Z)
\end{equation}
($\pi_1(Pic^{4g-4}(M))\cong \pi_1(Jac(M)) \cong
H_1(M,\Z)$, and $Pic_g^{4g-4}$ is a bundle over a contractible space,
so $\pi_1(Pic_g^{4g-4}) \cong H_1(M,\Z)$ as well).  An element of $\pi_1(Sym_g)$ will be a set of closed paths on a topological surface of genus $g$; $AJ_*$ will take the sum of these paths to their corresponding homology class. The composition of these maps must be trivial because the image of $\Q_g$ under the maps of (\ref{eq:seq}) is simply connected. 

Finally, consider the inclusion $\Q_\l \hookrightarrow \Q_g$, and the induced map of their fundamental groups. As previously noted it is possible that $\Q_\l$ is not connected, but for a particular connected component, $\Q_\l^0$ of $\Q_\l$, this gives us: 

\begin{equation}
\pi_1(\Q_\l^0) \to \pi_1(\Q_g) \to \pi_1(Sym_g) \to H_1(M, \Z).
\end{equation}

Since the first map is induced by an inclusion, we will refer to the composition of the first two maps as $i_*$ as well. Then we have:

\begin{prop} \label{prp:triv}
For all $\l \in \Lambda_g$ and any connected component $\Q_\l^0$ of $\Q_\l$
$$AJ_* \circ i_*:\pi_1(\Q_\l^0) \to H_1(M, \Z)$$  is trivial. 
\end{prop}

Note that the image of $\Q_\l \to Sym_g$ will be contained in
$Sym_g^\l$; thus, Proposition \ref{prp:triv} implies that the image of
$\pi_1(\Q_\l^0) \to \pi_1(Sym^\l_g)$ will be in the kernel of
$AJ_*$.  Since $\T_g$ is contractible, $\pi_1(Sym_g^\l) \cong B_\l$. Thus, we now use the generators of $B_\l$ constructed in the previous section to create a set of generators for $ker(AJ_*:\pi_1(Sym_g^\l) \to H_1(M, \Z))$. In subsection \ref{sub:nlg} we consider $\l=k^n$, and in subsection \ref{sub:gen} more general $\l$. 

\subsection{Strata with zeroes of only one weight} \label{sub:nlg}
We wish to calculate the kernel of $AJ_*$ for $\l=(k^n)$ and $n$ reasonably large. Recall that we defined $M_f$ to be the surface $M$ with $f$ punctures. We have the following from \cite{C}:

\begin{theorem} [Copeland] \label{thm:copeland}
If $M$ is a polyhedron (a two-dimensional cell complex) of genus $g$ with $n$ vertices and $f$ faces such that the associated graph has no double edges and no loops (edges with both ends at the same vertex), then $ker(AJ_*:\pi_1(Sym^n(M_f)) \to H_1(M_f, \Z))$ is generated by the edge set. Specifically the base point of $Sym^n(M_f)$ may be chosen to be the vertices of the cell complex, each face may be viewed as having a puncture, and each edge may be viewed as a transposition of its two vertices. 
\end{theorem}

 Theorem \ref{thm:copeland} implies that if we can construct a graph with any number of faces on $M$ of the form described in the theorem, then the kernel of $AJ_*:\pi_1(Sym^n(M)) \to H_1(M, \Z))$ is generated by transpositions.

Copeland shows in \cite{C} that it is possible to construct such a
graph for any $g>2$ and $n=4g-4$. We would like to show the same is
true for smaller $n$, since  $\l=(k^n)$ implies $n \le 4g-4$. In general, the best bound for $n$ we can hope to achieve will be on the order of $\sqrt{g}$. This is because a graph with $n$ vertices, no double edges, and no loops can have a maximum of $ \binom{n}{2}$ edges. The Euler characteristic then implies that the number of faces, $f$, of such a graph is given by:  
$$ 2-2g = n - \binom{n}{2} + f.$$
Since $f$ must be $\ge 1$, at best $g$ grows at the rate of $n^2$. More specifically, solving the equation above we get $n \ge \frac{3 + \sqrt{9 + 8(2g+f-2)}}{2}$. To show that graphs of the required form exist for $n$ close to this bound we will need some standard results from graph theory. Through the remainder of the section we will assume all graphs are connected, with no double edges or loops. 

An \textit{embedding} of a graph, $G$, into a surface of genus g, $M_g$, is a
homeomorphism $\phi:G \to M_g$. A \textit{2-cell embedding} of $G$ into $M_g$
is an embedding such that each component of $M_g \bs \phi(G)$ is a
2-cell. The \textit{genus} of $G$, $\g(G)$, is the minimal $g$ for
which $G$ embeds into $M_g$ (such an embedding will always be a 2-cell
embedding). The \textit{maximal genus} of $G$, $\g_M(G)$, is the the
maximal $g$ for which $G$ has a 2-cell embedding into $M_g$.  Finally,
let $K_n$ denote the complete graph on $n$ vertices, for $x \in \R$ let $\lceil x \rceil$ denote the smallest integer $\ge x$, and let $v(G), e(G),$ and $f(G)$ denote the number of vertices, edges, and faces of a graph. The following are well-known results in graph theory: 

\begin{prop} \label{prp:complete}
For any $n \in \mathbb{N}$:
\begin{enumerate}
\item  $\g_M(K_n)=\lceil \frac{e(K_n)-v(K_n)+1}{2} \rceil = \lceil \frac{(n-2)(n-1)}{4} \rceil$.
\item $\g(K_n)=\lceil \frac{(n-3)(n-4)}{12} \rceil$. 
\item $K_n$ has a 2-cell embedding into $M_g$ if and only $\g(K_n) \le g \le \g_M(K_n)$. 
\end{enumerate}
\end{prop}

 See \cite{KRW}, for example, for a survey of these results. Note that
 the maximal genus of $G$ is simply the largest genus for which $v(G) - \binom{v(G)}{2} + f = 2 - 2g$ has a positive solution for $f$. A particular 2-cell embedding of $K_n$ into $M_g$ has $\binom{n}{2} - n + 2 -2g$ faces. We would like to show that graphs with a wider range of faces embed into $M_g$, and in fact knowing that $K_n$ 2-cell embeds into $M_g$ we can also show that `almost complete' graphs on $n$ vertices have 2-cell embeddings into $M_g$. 

\begin{lemma} \label{lem:face}
If $K_n$ has a 2-cell embedding into $M_g$, then there also exists a graph with $n$ vertices and any positive number of faces $\le \binom{n}{2} - n + 2 -2g$ that has a 2-cell embedding into $M_g$. 
\end{lemma}

\begin{proof}
We prove the lemma by induction on the number of faces. 
First, $K_n$ has $\binom{n}{2}-n+2-2g$ faces and embeds into $M_g$. Now suppose we have a connected graph $G$ with $f(G)$ faces that has a 2-cell embedding $\phi:G \to M_g$, and suppose that $2 \le f(G) \le \binom{n}{2}-n+2-2g$. Each edge of a graph is adjacent to two faces (it may be adjacent to the same face twice), and since $G$ is connected any face must share at least one edge with some other face. Call this edge $e$. $G \bs e$ will still be connected and $\phi|_{G \bs e}$ will be a 2-cell embedding of $G \bs e$ with one face fewer than $G$. 
\end{proof}

\begin{prop} \label{prp:graph}
 Let $1 \le f \le 4g-4$. Then for any $n \ge \frac{3 +
   \sqrt{9 + 8(2g+f-2)}}{2}$ there exists a graph on a surface of
 genus $g$ with $n$ vertices, $f$ faces, no loops and no double edges. 
\end{prop}

\begin{proof}
First, it suffices to show that for $n = \lceil \frac{3 + \sqrt{9 + 8(2g+f-2)}}{2} \rceil$ there exists such a graph, $G$; to construct such a graph for $n' > n$ we simply subdivide the edges of $G$ with the required number of additional vertices.

 Let $e$ be the number of edges of $G$. Notice that $n = \frac{3 + \sqrt{9 + 8(2g+f-2)}}{2}$ implies $e=\binom{n}{2}$, while $n-1= \frac{3 + \sqrt{9 + 8(2g+f-2)}}{2}$ implies $e=\binom{n-1}{2}+1$. Thus $n=\lceil \frac{3 + \sqrt{9 + 8(2g+f-2)}}{2} \rceil$ implies $\binom{n-1}{2}+1 < e \le \binom{n}{2}$. By assumption $1 \le f \le 4g-4$, and since $n-e+f=2-2g$, $g = \frac{e-n-f+2}{2}$. Thus: 
$$\frac{\binom{n-1}{2}+1 - n - (4g-4) + 2}{2} \le g \le \frac{\binom{n}{2} - n -1 +2}{2}.$$
Simplifying:
$$\frac{n^2-5n+14}{12} \le g \le \frac{(n-1)(n-2)}{4}$$ 
Then by Proposition \ref{prp:complete}, for $n = \lceil \frac{3 + \sqrt{9 + 8(2g+f-2)}}{2} \rceil$ and $1 \le f \le 4g-4$, $K_n$ has a 2-cell embedding into $M_g$. If  $n = \frac{3 + \sqrt{9 + 8(2g+f-2)}}{2}$, $K_n$ has $f$ faces and is the desired graph. Otherwise $K_n$ has more than $f$ faces, but by Lemma \ref{lem:face} a graph with $n$ vertices and any number of faces fewer than $f(K_n)$ is also embeddable in $M_g$. This proves the proposition.  
\end{proof}

\begin{corollary} \label{cor:imv2}
For $g \ge 2$, $1\le f \le 4g-4$ and $n \ge \frac{3 + \sqrt{9 +8(2g+f-2)}}{2}$, $ker(AJ_*:\pi_1(Sym^n(M_f)) \to H_1(M_f, \Z))$ is generated by transpositions.
\end{corollary}

\begin{proof}
This is an immediate consequence of Theorem \ref{thm:copeland} and Proposition \ref{prp:graph}
\end{proof}

\begin{corollary} \label{cor:image}
For $g \ge 2$, $n \ge \frac{3 + \sqrt{1+16g}}{2}$, and $\l=(k^n)$, $ker(AJ_*:\pi_1(Sym_g^\l) \to H_1(M, \Z))$ is generated by transpositions.
\end{corollary}

\begin{proof}
Substitute $1$ for $f$ in the formula of Corollary \ref{cor:imv2} and note that $ker(AJ_*:\pi_1(Sym^\l(M)) \to H_1(M, \Z))$ is a subgroup of the kernel of $AJ_*:\pi_1(Sym^n(M_1)) \to H_1(M_1, \Z)$.
\end{proof}

Corollary \ref{cor:image} gives us the structure of the kernel of $AJ_*$ for $\l=(k^n)$ and $n$ reasonably large. 

For $\l$ with only a few points, it is more difficult to enumerate a set of generators for the kernel of $AJ_*$. It is not true in general that the kernel will be generated by transpositions; for example, consider the stratum $\Q_g(4g-4)$. We may move the single marked point around a curve that is homologically but not homotopically trivial. This will be in the kernel of $AJ_*$ but is not a product of transpositions. However, for $\l$ with sufficiently many zeroes of the same order we may make some generalizations to the results of this subsection. In particular, a combination of Corollary \ref{cor:imv2} and the Fadell-Neuwirth fibration will allow us to generalize Corollary \ref{cor:image} to a larger class of $\l$, and this is what we will do in the next subsection. 

\subsection{Strata with zeroes of more than one weight} \label{sub:gen}

We would like to generalize Corollary \ref{cor:image} to $\l=(k_1^{n_1},...,k_l^{n_l})$ with $n_1$ large. From Section \ref{sec:pi1} we know that for general $\l$, $B_\l \cong \pi_1(Sym_g^\l)$ is generated by transpositions or square transpositions of $p_i$ with $p_j$, $\sigma_{ij}$ or $\kappa_{ij}$, and moving a point $p_i$ around $l_r \in \pi_1(M)$, $\rho_{ir}$. We can immediately show that some of these generators are in $ker(AJ_*:\pi_1(Sym_g^\l) \to H_1(M, \Z))$:

\begin{lemma} \label{lem:trp} 
Any transposition or square transposition of points in $\pi_1(Sym_g^\l)$ is in the kernel of $AJ_*$.
\end{lemma}

\begin{proof}
A transposition of two points of equal weight consists of moving them in opposite directions along homotopic paths. The sum of these paths is then homotopic (and therefore homologous) to zero. A square transposition of two points of unequal weight moves each point along some path, and then back along a homologous path. Thus both points follow paths that are homologous to zero.   
\end{proof}

An individual $\rho_{ir}$ will not be in the kernel of $AJ_*$; however, there are two cases when it is easy to see that a product of them will be. First, if two sets of points of equal total weight follow $l_r$ and $l_r^{-1}$ respectively, then their paths will cancel each other out in $H_1(M, \Z)$. Second, a single $p_i$ may follow a path that is homologically trivial but not homotopically trivial. More precisely, we have the following two definitions. 

\begin{definition} Fix $r$, $1 \le r \le 2g$. A \textit{null $\rho_r$} is a
  product $\Pi_{n=1}^m \rho_{i_n r}^{\pm 1}$ such that $\sum_{n=1}^m \pm k_{i_n}$ $= 0$, where the sign in front of $k_{i_n}$ is given by the sign of the exponent of $\rho_{i_n r}$.
\end{definition}

\begin{definition} Let $\l=(k_1,...,k_n)$, fix $i$, $1 \le i \le n$, and let $M_{n-i}$ be $M$ punctured at the zeroes of $q$ of weight $k_{i+1},...,k_n$. An \textit{$i$-commutator} is a product of $\rho_{ij}$ and $\kappa_{il}$, $1 \le j \le 2g$, $i+1 \le l \le n$, such that the path followed by $p_i$ is in $[\pi_1(M_{n-i}), \pi_1(M_{n-i})]$. 
\end{definition}
Any null $\rho_r$ or $i$-commutator is in $ker(AJ_*)$.  We show that for $\l=(k_1^{n_1},...,k_l^{n_m})$ with $n_1$ large, transpositions, square transpositions, null $\rho_r$, and $i$-commutators for $p_i$ not of order $k_1$ suffice to generate the kernel of $AJ_*$. For $n_2,...,n_m$ all sufficiently large, transpositions, square transpositions, and null $\rho_r$ will suffice to generate.

To show this first recall that for any $\l$ of length $n$: 
 
\begin{equation} \label{eq:ses1}
SB_n \stackrel{\iota}{\to} \pi_1(Sym_g^\l) \cong B_\l \stackrel{pr}{\to} S_\l
\end{equation}
is a short exact sequence. From Section \ref{sec:pi1} we have that $SB_n$ is generated by $\rho_{ir}$ and $\kappa_{ij}$, $1 \le i<j \le n$, $1 \le r \le 2g$. Let $K=ker(AJ_*)$ and let $K'=ker(AJ_* \circ \iota:SB_{n} \to H_1(M, \Z))$. 

\begin{lemma} \label{lem:tech0}
Any $Z \in \pi_1(Sym_g^\l)$ can be written as $Y \cdot X$ where $Y$ is a product of transpositions and $X \in im(\iota:SB_n \to \pi_1(Sym_g^\l))$. 
\end{lemma}

\begin{proof}
Let $pr(Z) = \bar{Y} \in S_\l$. For every $p_i, p_j$ of the same weight, pick a transposition, $\sigma_{ij} \in \pi_1(Sym_g^\l)$. Let $pr(\sigma_{ij})=\bar{\sigma}_{ij}$ and note that $S_\l$ is generated by the $\bar{\sigma}_{ij}$. Then we can write $\bar{Y}$ as a product of the $\bar{\sigma}_{ij}$, and we construct $Y \in \pi_1(Sym_g^\l)$ by writing $Y$ as a product of the corresponding $\sigma_{ij}$. Then $pr(Y)=pr(Z)=\bar{Y}$. Let $X = Y^{-1}Z$. Since $pr(X)=1$ and (\ref{eq:ses1}) is exact, $X \in im(SB_n \to \pi_1(Sym_g^\l))$.
\end{proof}

Thus to prove that any $Z \in K$ is a product of null $\rho_r, \sigma_{ij}$, $\kappa_{ij}$, and $i$-commutators, it suffices to show that this is true for all $X=\iota(X')$, $X' \in K'$. Since $K' < SB_n$ we may use the short exact sequence of Theorem \ref{thm:fn}: 

\begin{equation} \label{eq:ses2}
SB_{a, b} \stackrel{\iota'}{\to} SB_{a+b}=SB_n \stackrel{pr'}{\to} SB_{b}
\end{equation} 
for any $a,b$ such that $a+b=n$. Let $K''=ker(AJ_* \circ \iota \circ \iota':SB_{a,b} \to H_1(M, \Z))$. Lemma \ref{lem:tech3} will prove that words in $K''$ can be written as a product of the desired elements, and Theorem \ref{thm:gen} will then prove the same for words in $K'$ (and therefore $K$). 

First we need a technical lemma, and in it we break up $SB_{a,b}$ further to analyze it. Again apply Theorem \ref{thm:fn} to get:
\begin{equation} \label{eq:ses3}
SB_{a-1,b+1} \stackrel{\iota''}{\to} SB_{a,b} \stackrel{pr''}{\to} SB_{1,b}
\end{equation}

\begin{lemma} \label{lem:tech1} Let $\l=(k_1^{n_1},....,k_{m}^{n_m})$, with $\sum n_i = n$. Pick any $a$, $1 \le a \le n$, let $b=n-a$, and let $p_a$ have weight $k_l$, $1 \le l \le m$. Then $im(pr'':K''<SB_{a,b} \to SB_{1,b})$ lies in $\{S \in SB_{1,b} | AJ_*(S) \in d H_1(M, \Z) \}$, where $d$ is the smallest positive integer such that there exist $c_1, c_2,...,\hat{c_l},...,c_m \in \Z$ such that $c_1 \cdot k_1 + c_2 \cdot k_2 +...+d \cdot k_l+...+c_{m} \cdot k_{m} = 0$ has a solution. 
\end{lemma}

\begin{proof}
$SB_{a,b}$ is generated by $\rho_{ir}, \kappa_{ij}$, $1 \le i \le a$, $1 \le j \le a+b$, $1 \le r \le 2g$. We need not consider the $\kappa_{ij}$ as they all go to zero under $AJ_*$. For a fixed $r$ the product of $\rho_{ir}^{\pm 1}$ in $K''$ must be such that $\sum \pm k_i=0$. This implies that the number of times $\rho_{ar}$ occurs in a particular word must be a multiple of $d$. 
\end{proof}

\begin{lemma} \label{lem:tech3}
Let $\l$, $c_1, c_2,...,\hat{c_l},...,c_m$, and $d$ be as in Lemma \ref{lem:tech1}, with $p_a$ again of weight $k_l$. For all $X'' \in K''$ there exists $W \in K''$ such that $pr''(X'') = pr''(W)$ and $\iota \circ \iota'(W) \in \pi_1(Sym_g^\l)$ is a product of null $\rho_r$, $\kappa_{aj}$ or $\sigma_{aj}$, and $a$-commutators.  
\end{lemma}

\begin{proof}
$SB_{1,b}=\pi_1(M_b)$ is generated by $\bar{\rho}_{ar},\bar{\kappa}_{aj}$, $a+1 \le j \le a+b$, $1 \le r \le 2g$, where  $\rho_{ar},\kappa_{aj} \in SB_{a,b}$ project to their corresponding barred elements. (In $SB_{a,b}$ the points $p_1,...,p_a$ move around punctures $p_{a+1},...,p_{a+b}$; under $pr''$, $p_a$ moves around $p_{a+1},...,p_{a+b}$.)  Let $G=[SB_{1,b}, SB_{1,b}]$ be the commutator of $SB_{1,b}$. 

Since $G$ abelianizes $SB_{1,b}$, there exists $h \in G$ such that $h \cdot pr''(X'')$ is a product of $\bar{\rho}_{ar}, \bar{\kappa}_{aj}$, such that all $\bar{\kappa}_{aj}$ are on the right of all $\bar{\rho}_{ar}$ and if $r_1 < r_2$, $\bar{\rho}_{a r_1}$ is to the left of $\bar{\rho}_{a r_2}$. In other words, $h \cdot pr''(X'')$ is a word such that for fixed $r$ all $\bar{\rho}_{ar}$ are adjacent. By Lemma \ref{lem:tech1}, the power of any $\bar{\rho}_{ar}$ in $h \cdot pr''(X)$ must be a multiple of $d$. 

Now we construct $W \in K''$ by inserting $c_i$ elements of weight $k_i$ moving around $l_r$ adjacent to each set of $d$ elements of weight $k_l$ moving around the $l_r$ in $h \cdot pr''(X'')$, and then multiplying by $h^{-1}$. For example, let $\l=(k_1, k_2)$ with $p_1$ of weight $k_1$, $p_2$ of weight $k_2$, and $c k_1 - d k_2 = 0$. If $h \cdot pr''(X)= \bar{\rho}_{1r_1}^{n_1} \bar{\rho}_{1r_2}^{n_2} \bar{\kappa}_{12}$, then $W = h^{-1}(\rho_{2r_1}^{n_1} \rho_{1r_1}^{-\frac{c}{d} n_1}) (\rho_{2r_2}^{n_2} \rho_{1r_2}^{-\frac{c}{d} n_2})\kappa_{12}$, where Lemma \ref{lem:tech1} implies $n_1, n_2$ are divisible by $d$. Both of the elements in parentheses are null $\rho_r$. $W$ is thus made up of null $\rho_r$, $\kappa_{aj}$, and $a$-commutators. Taking the inclusion of $W$ in $K$ under $\iota \circ \iota'$ does not change this, and by construction $pr''(W)=pr''(X'')$.
\end{proof}

Notice that if there are sufficiently many points of weight $k_l$ (greater than or equal to $\frac{3 + \sqrt{9 + 8(2g+b-2)}}{2}$), then $a$-commutators can be written as a product of transpositions by Corollary \ref{cor:imv2}. 

\begin{theorem} \label{thm:gen}
Let $\l=(k_1^a,k_2^{b_2},...,k_m^{b_m})$, with $\sum_2^m b_i=b$ and $a \ge \frac{3 + \sqrt{9 +8(2g+b-2)}}{2}$ (as in Corollary \ref{cor:imv2}). Then any $Z \in K$ can be written as a product of null $\rho_r$, $\sigma_{ij}$, $\kappa_{ij}$, and $i$-commutators for $p_i$ not of weight $k_1$. 
\end{theorem} 

\begin{proof}
By Lemma \ref{lem:tech0} it suffices to prove the theorem for any $X=\iota(X')$ such that $X' \in K' < SB_{a+b}$. By Corollary \ref{cor:imv2} the theorem is true for for $\iota \circ \iota'(X'')<B_{a,b}$, $X'' \in K'' < SB_{a,b}$. To get from $SB_{a,b}$ to $SB_{a+b}$ we use Lemma \ref{lem:tech3} as follows. Starting with any $X'=X_0 \in SB_{a+b}$ we construct a sequence of $X_k \in SB_{a+b-k,k}$, $0 \le k \le b$, $X_b \in SB_{a,b}$, such that if $X_{k+1}$ is generated by elements of the desired form then $X_k$ is as well. Since $X_b$ is generated by elements of the desired form, this will imply that $X_0$ is as well, thus proving the theorem.  

In particular, let $\iota_b:SB_{a,b} \stackrel{\iota'}{\to} SB_{a+b} \stackrel{\iota}{\to} \pi_1(Sym_g^\l)$ be the composition of the two injections given in (\ref{eq:ses1}), (\ref{eq:ses2}). Let $K_b=K''=ker(AJ_* \circ \iota_b:SB_{a,b} \to H_1(M, \Z))$, and pick $X_b \in K_b$. Then $\iota_b(X_b) \in (K \cap B_{a,b}) \subset \pi_1(Sym_g^\l(M))$. By Corollary \ref{cor:imv2} $i_b(X_b)$ can be written as a product of transpositions. Similarly for any $k$, $0 \le k < b$, define $\iota_k:SB_{a+b-k,k} \to SB_{a+b} \to \pi_1(Sym_g^\l)$. Also define $K_k=ker(AJ_* \circ \iota_k:SB_{a+b-k,k} \to H_1(M, \Z))$. 

Fix $k_0$, $0 \le k_0 < b$, and assume any element in $\iota_{k_0+1}(K_{k_0+1})$ can be written as in the statement of the theorem.  Then consider the following version of the short exact sequence in (\ref{eq:ses3}): 
$$SB_{a+b-(k_0+1),k_0+1} \stackrel{\iota''}{\to} SB_{a+b-k_0,k_0} \stackrel{pr''}{\to} SB_{1, k_0}$$
Let $X_{k_0}$ be an element of $K_{k_0}$. By Lemma \ref{lem:tech3} we construct $W_{k_0} \in K_{k_0}$ such that $pr''(W_{k_0})=pr''(X_{k_0})$ and $\iota_{k_0}(W_{k_0})$ is a product of null $\rho_r$, $\kappa_{(a+b-k_0)j}$ or  $\sigma_{(a+b-k_0)j}$ and $(a+b-k_0)$-commutators ($p_{(a+b-k_0)}$ is the `single point' of $SB_{1,k_0}$).  

Now define $\bar{X}_{k_0+1}=W_{k_0}^{-1}X_{k_0}$ and note that
$pr''(\bar{X}_{k_0+1})=1$. Thus $\bar{X}_{k_0+1} \in im(i'':SB_{a+b-(k_0+1),k_0+1}
\to SB_{a+b-k_0,k_0})$. Since $\iota''$ is injective we may define  $X_{k_0+1} \in SB_{a+b-(k_0+1),k_0+1}$ to be $\iota''^{-1}(\bar{X}_{k_0+1})$. Since $\iota''$ is just an inclusion,
$\iota_{k_0}(\bar{X}_{k_0+1})$ is equal to $\iota_{k_0+1}(X_{k_0+1})$ which by assumption may be written as a product of the elements in the statement of the theorem. 

Consequently $\iota_{k_0}(X_{k_0})$ can be written as a product of $\iota_{k_0}(W_{k_0})$ and $\iota_{k_0}(\bar{X}_{k_0+1})$ both of which may be written as products of elements in the statement of the theorem. 

In particular, for $k_0=0$ we get that for any $X_0=X' \in K_0 = K'$, $\iota(X')$ may be written as in the statement of the theorem.  
\end{proof}

\begin{corollary} \label{cor:gen2}
Let $\l=(k_1^a,k_2^{b_2},...k_n^{b_m})$ with $\sum_2^m b_i=b$ and $a \ge \frac{3 + \sqrt{9 + 8(2g+b-2)}}{2}$.  Additionally assume that $b_i \ge \frac{3 + \sqrt{9 + 8(2g+b_{i+1}+...+b_{m}-2)}}{2}$, $2 \le i \le m$. Then any $Z \in ker(\pi_1(Sym_g^\l) \to H_1(M, \Z))$ may be written as a product of $\sigma_{ij}$, $\kappa_{kl}$, and  null $\rho_r$.
\end{corollary}

\begin{proof}
This differs from Theorem \ref{thm:gen} only in that we eliminate the $k$-commutators, $a < k \le a+b$. We may do this because all of these elements are contained in $ker(B_{b_i,b_{i+1}+...+b_m} \to H_1(M, \Z))$ and may thus by Corollary \ref{cor:imv2} be written as products of transpositions.  
\end{proof}

This gives us a description of $ker(AJ_*:Sym_g^\l \to H_1(M, \Z))$ for $\l$ with at least one reasonably large set of points of the same weight (on the order of $\sqrt{g}$). The rest of the paper will be spent in determining which of the elements of this kernel are also elements of $\pi_1(\Q_\l)$. To do so we first present some methods for surgering existing quadratic differentials to create new ones.  
 
\section{Some Local Surgeries} \label{sec:adj}
We denote by $\Q_\l^0$ a connected component of $\Q_\l$, and by $cl(\Q_\l)$ the closure of $\Q_\l$. 
Choose $\l_1, \l_2$ such that $\Q_{\l_2}^0 \subset cl(\Q_{\l_1}^0)$. In this section we present a method whereby, for certain $(M, q) \in \Q_{\l_1}^0$ sufficiently near to an element of $\Q_{\l_2}^0$, we may construct some transpositions and square transpositions of zeroes of $q$ in $\pi_1(\Q_\l^0,(M,q))$. For $\l=(k_1,...,k_n)$ we refer to the zeroes of $q$ as $p_1,...,p_n$ with $p_i$ of order $k_i$, and we let $P=\{p_1,..,p_n\}$. We also sometimes refer to a (square) transposition of $p_i$ and $p_j$ - by this we will mean a transposition if $k_i=k_j$ and a square transposition otherwise.

First we consider the following two lemmas, which closely follow lemmas in \cite{EMZ} and \cite{L1}, amongst others, and which allow us in certain cases to construct an element of $\Q_{\l_1}^0$ from an element of $\Q_{\l_2}^0$.
\begin{lemma} \label{lem:two} 
Let $(M,q) \in \Q_g(k_1,...,k_n)=\Q_{\l_2}$. Pick $l_1, l_2$ such that
$l_1+l_2=k_i$, where if $k_i$ is even then $l_1$ and $l_2$ are even as well. Then it is possible to construct a deformation of $(M,q)$, $(M',q') \in
\Q_{\l_1}=\Q_g(k_1,...,k_{i-1},l_1,l_2,k_{i+1},...,k_n)$, such that the flat
metric on $M$ is unchanged outside of an $\epsilon$ disk around $p_i$. 
\end{lemma}

\begin{figure}
\begin{center}
\includegraphics[width=.3\textwidth]{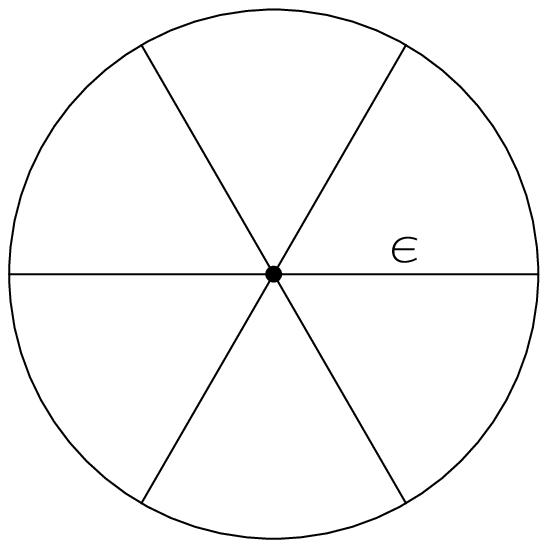}
\ \ \ \ \ \
\includegraphics[width=.3\textwidth]{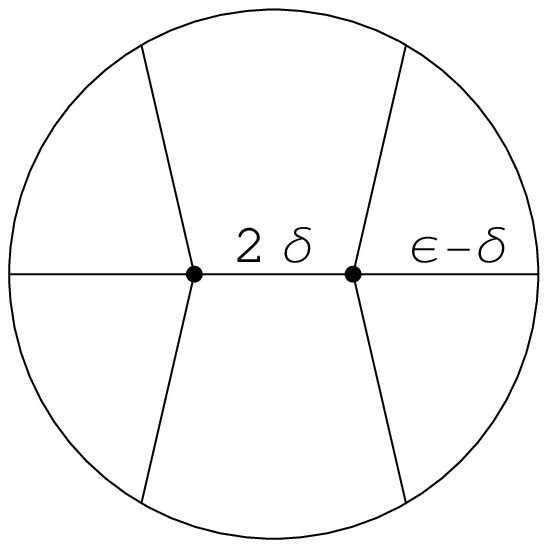}
\\
\includegraphics[width=.3\textwidth]{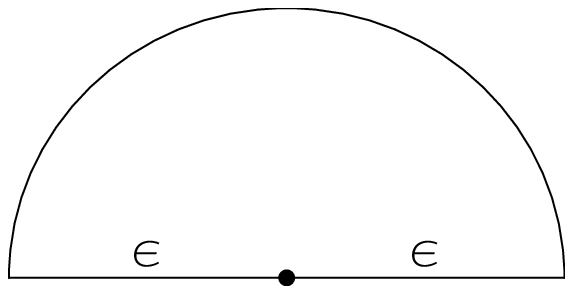}
\ \  
\includegraphics[width=.3\textwidth]{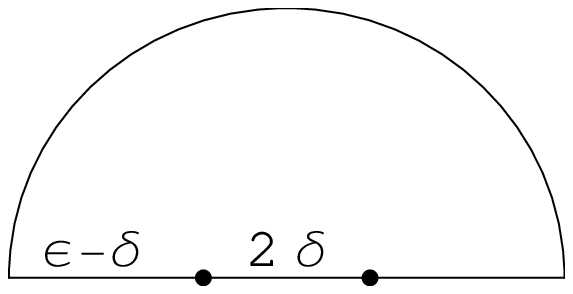}
\ \ 
\includegraphics[width=.3\textwidth]{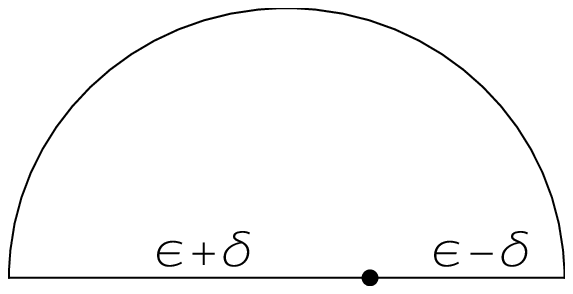}

\end{center}
\caption{Taking a zero of order 4 to two zeroes of order 2 via
  splitting/moving the zero on individual half-disks. These diagrams are purely schematic and not to scale; the angle between any two straight lines in the disks is $\pi$.}
\label{fig:spliteven}
\end{figure}

\begin{lemma} \label{lem:three} Let $(M,q) \in \Q_g(k_1,...,k_n)$. Pick $l_1,l_2,l_3$ of any parity, such that $l_1+l_2+l_3=k_i$.  Then it is possible to construct a deformation of $(M,q)$, $(M',q') \in \Q_g(k_1,...,
k_{i-1},l_1,l_2,l_3,k_{i+1},...,k_n)$, such that the flat metric on $M$ is unchanged outside of an $\epsilon$ ball around $p_i$. Similarly, given any $l_1,l_2,l_3,l_4$ such that $l_1+l_2+l_3+l_4=k_i$ then there exists $(M',q')  \in \Q_g(k_1,...,k_{i-1},l_1,l_2,l_3,l_4,k_{i+1},...,k_n)$ with flat metric unchanged outside of a ball around $\ep$. 
\end{lemma} 

\begin{figure} 
\begin{center}
\includegraphics[width=.4\textwidth]{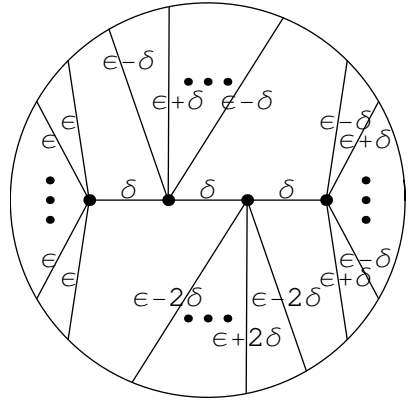}
\end{center}
\caption{Splitting an even zero of high order into 4 odd zeroes of lower order.}
\label{fig:four}
\end{figure}

\begin{proof}[Sketch of Proof]
A proof of all but the 4 point case may be found in \cite{L1}; we prove the two point case with both zeroes even here. Consider a disk $D_\epsilon$, centered at $p_i$, of radius $\epsilon$ in the flat metric given by $q$, where $\ep$ is small enough that $D_\epsilon$ contains no other zeroes of $q$. Since $p_i$ is of order $k_i$, $k_i+2$ horizontal trajectories will dead-end into $p_i$. Cut along these to make $k_i+2$ half-disks, each with a marked point given by $p_i$ half way along the
cut. Pick some $\delta$ with $0<\delta<\epsilon$. If $l_1$ and $l_2$ are both even, construct two special half-disks by splitting the marked point given by $p_i$ into two marked points $2\delta$ apart. Also shift the marked point given by $p_i$ by $\delta$ on the remaining disks, as in Figure \ref{fig:spliteven}. The trajectory structures on the individual half-disks do not change, and the half-disks can be glued back together, again as in Figure \ref{fig:spliteven}, to give an $\epsilon$ disk containing zeroes of order $l_1,l_2$, with the same trajectory structure at the boundary as the original $D_\epsilon$. The remaining two and three point cases are proved similarly. 

For the case of four zeroes, if any of the $l_j$ are even then we may apply some combination of two and three point surgeries to get the desired result. If $l_1,l_2,l_3,l_4$ are all odd then we cut $D_\ep$ around $p_i$ into $k_i+2$ half-disks and glue as in Figure \ref{fig:four}. Again, in all of these surgeries we have a choice of both $\delta$ and $\theta$. 
\end{proof}

Define $\Q_{\l_1}^0$ to be \textit{adjacent} to $\Q_{\l_2}^0$ if it is possible to obtain an element in $\Q_{\l_1}^0$ by (possibly repeatedly) applying the surgeries of Lemmas \ref{lem:two} and \ref{lem:three} to an element in $\Q_{\l_2}^0$. This puts a poset ordering on $\Lambda$, the set of partitions of $4g-4$: we say $\l_1 > \l_2$ if there exists a component of $\Q_{\l_1}$ that is adjacent to a component of $\Q_
{\l_2}$. Note that if $\Q_{\l_1}^0$ is adjacent to $\Q_{\l_2}^0$, then $Q_{\l_2}^0 \subset cl(\Q_{\l_1}^0)$. 

Lemmas \ref{lem:two} and \ref{lem:three} detail surgeries that allow us to break up zeroes of quadratic differentials; however, again following \cite{EMZ} we may \textit{collapse} zeroes in the reverse process. We say three points, $p_1, p_2$, and $p_3$, are \textit{co-linear} if for a fixed $\theta$ there exist $\theta$-trajectories between $p_1,p_2$ and $p_2, p_3$. 

\begin{lemma}\label{lem:col}
Suppose $(M',q')\in \Q_g(k_1,...,k_n)$ is such that there exists a saddle connection between $p_1$ and $p_2$ of length $\delta$, all other saddle connections from $p_1$ or $p_2$ are of length at least $3 \delta$, and $k_1,k_2$ are not both odd. Then there exists $(M,q) \in \Q_g(k_1+k_2,k_3,...,k_n)$ such that $(M',q')$ may be obtained from $(M,q)$ via the surgery of Lemma \ref{lem:two}. 

Suppose $(M',q')\in \Q_g(k_1,...,k_n)$ is such that there exist saddle connections between $p_1, p_2$ and $p_2, p_3$ of length $\delta$, $p_1, p_2, and p_3$ are co-linear, and all saddle connections from $p_2$ are of length at least $4 \delta$. Then there exists $(M,q) \in \Q_g(k_1+k_2+k_3,...,k_n)$ such that $(M',q')$ may be obtained from $(M,q)$ via the surgery of Lemma \ref{lem:two}. 
\end{lemma}

\begin{figure}
\begin{center}
\includegraphics[width=.3\textwidth]{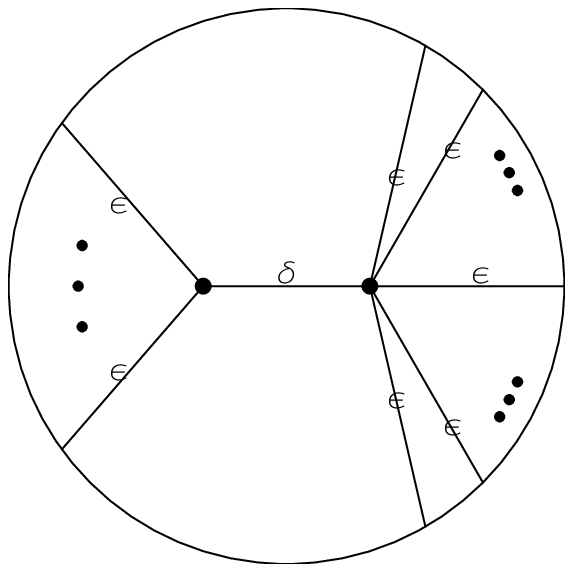}
\ \ \ \
\includegraphics[width=.3\textwidth]{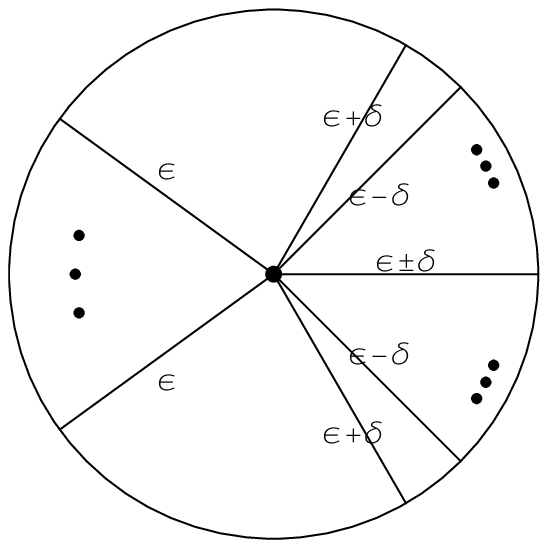}
\\
\includegraphics[width=.3\textwidth]{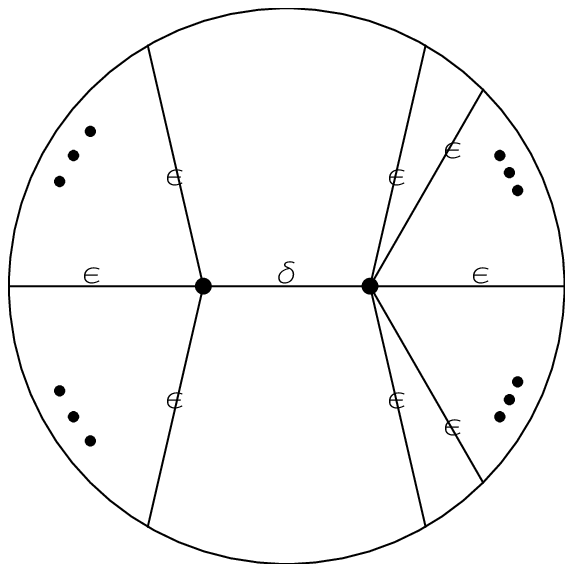}
\ \ \ \
\includegraphics[width=.3\textwidth]{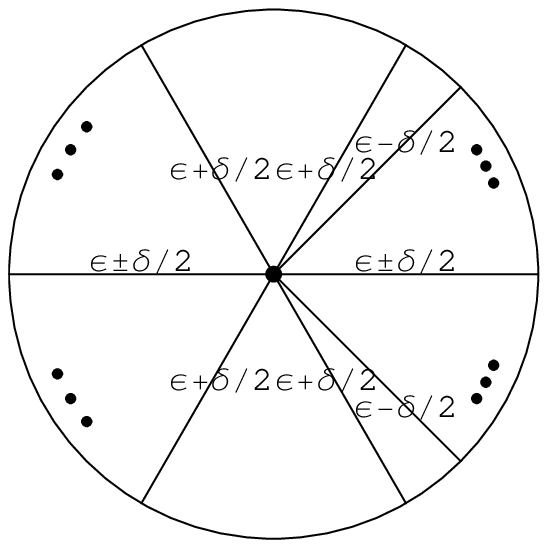}
\\
\includegraphics[width=.3\textwidth]{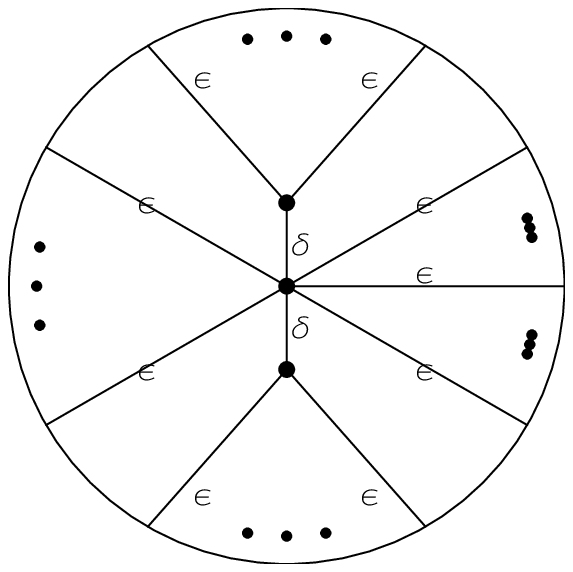}
\ \ \ \ 
\includegraphics[width=.3\textwidth]{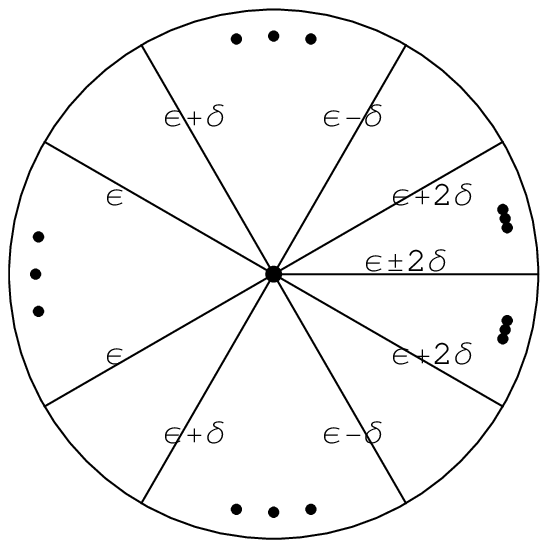}
\\
\includegraphics[width=.3\textwidth]{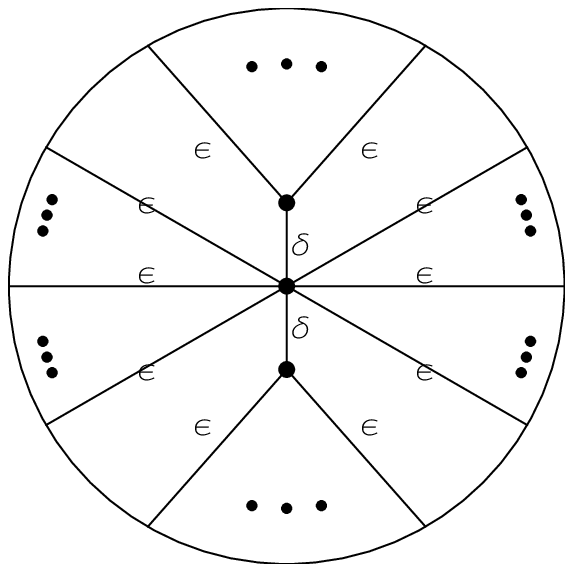}
\ \ \ \
\includegraphics[width=.3\textwidth]{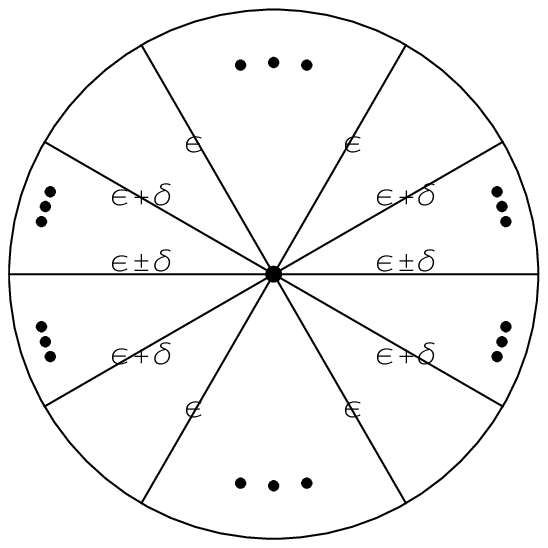}
\\
\end{center}
\caption{i. Colliding two zeroes, one of odd order and one of even order. ii. Colliding two zeroes of even order. iii. Colliding three zeroes of odd order. iv. Colliding three zeroes, two of any order and one even.} 
\label{fig:coll}
\end{figure}

\begin{proof}
We first prove the two point case. Label the points in each of the two top left disks of Figure \ref{fig:coll} by $p_1$ and $p_2$. For $\ep = 2 \delta$ both of the two top left disks are contained in $D_{3\delta}(p_1) \cup D_{3\delta}(p_2)$ and therefore contain no other zeroes of $q'$. Since $\ep > \delta$ we may reverse the cutting and pasting process of Lemma \ref{lem:two} to obtain $(M,q)$. 
For the three point case, label the points in the the two bottom left discs of Figure \ref{fig:coll} by $p_1, p_2, p_3$, with $p_2$ in the middle. For $\ep =3 \delta$ $D_{4\delta}(p_2)$ contains either of these two discs.  Thus these disks contain no other zeroes of $q'$ and since $\ep > 2\delta$ we may reverse the cutting and pasting of Lemma \ref{lem:three} to obtain $(M,q)$.  
\end{proof}

In the next proposition we will construct (square) transpositions of zeroes of $(M', q') \in \Q_{\l_1}$ by colliding two zeroes to get $(M,q)$ in an adjacent $\Q_{\l_2}$ and then breaking up the newly formed zero of $q$ with respect to varying $\theta$. Breaking up a single zero into two zeroes in this manner gives us a transposition; however, breaking a single zero into three gives us something slightly different.  Define a \textit{transposition of $p_i$ and $p_j$ mod $p_k$} to be an element of $\pi_1(\Q_\l, (M,q))$ such that $p_i, p_j$ follow paths that are not homotopic on $M \bs P$, but are homotopic on $M \bs P \cup p_k$. If $e$ is a path between $p_i, p_j$ through $p_k$ we will refer to a transposition of $p_i, p_j$ mod $p_k$ along $e$; by this we mean that $p_i$ and $p_j$ follow small deformations of $e$ on opposite sides of $p_k$. Similarly if we
refer to a transposition of $p_i, p_j$ along $e$ we mean their transposition along a slight deformation of $e$ to one side of $p_k$.

\begin{prop} \label{prp:tra}Let $\l_1=(k_1,k_2,...,k_n)$, $\l_2=(k_1+k_2,k_3,...,k_n)$, and suppose $(M',q') \in \Q_{\l_1}$ satisfies the conditions of Lemma \ref{lem:col} in the two point case. Then there exists a trajectory of length $\delta$ between $p_1$ and $p_2$ (specified by \ref{lem:col}), and an element of $\pi_1(Q_{\l_1}, (M',q'))$ corresponding to the (square) transposition of $p_1, p_2$ along the specified trajectory.  

Suppose $\l_1=(k_1,k_2,...,k_n)$, $\l_3=(k_1+k_2+k_3,k_4,...,k_n)$, and $(M',q') \in \Q_{\l_1}$ satisfies the conditions of Lemma \ref{lem:col} in the three point case.  Then there exists a transposition of $p_1, p_2$ mod $p_3$ along the union of two trajectories specified by \ref{lem:col}.
\end{prop}
 
\begin{proof} We first prove the proposition for the two point case. Suppose $(M',q')$ is obtained from $(M,q) \in \Q_{\l_2}$ by surgering with respect to $\delta_0, \theta_0$. Let $e_0$ be the new trajectory between $p_1, p_2$ of length $\delta_0$. Create a curve $\eta:[0,\delta_0] \to Cl(\Q_{\l_1})$ by surgering $(M,q)$ with respect to $t, \theta_0$, $t \in (0, \delta_0]$.  Now, there exists a ball of some radius $\alpha$ around the zero of $q$ of order $k_1+k_2$ containing no other zeroes of $q$, so we apply the surgery of Lemma \ref{lem:two} to $(M, q)$ with respect to some fixed $\delta_1 < \textrm{min}\{\delta_0, \alpha\}$ and vary $\theta$ between $0$ and either $\pi$ or $2\pi$ (depending on whether we transpose or square transpose $p_1$ and $p_2$) to get a loop of surfaces, $\eta' \subset \Q_{\l_1}$, centered at $(M, q)$. $\eta$ intersects $\eta'$ at $\eta(\delta_1)$, and we define $\tilde{\eta}$ to be the sub-curve of $\eta$ from $\delta_0$ to $\delta_1$. Then $\tilde{\eta} \circ \eta' \circ \tilde{\eta}^{-1}$ gives us the desired loop in $\pi_1(\Q_{\l_1}, (M', q'))$. 

Following the same procedure for three points, we get a (square) transposition of $p_1, p_2$ mod $p_3$ along the union of the two trajectories between them in $\pi_1(Q_{\l_1}, (M',q'))$. 
\end{proof}
Proposition \ref{prp:tra} allows us to create (square) transpositions of
various zeroes of a quadratic differential (possibly mod a third zero); however, the (square) transpositions created this way may be based at different elements of $\Q_\l$. In the next section we show that one way of obtaining elements based at the same $(M,q)$ is to consider hyperelliptic quadratic differentials.  

\section{ (Square) transpositions in certain strata} \label{sec:hyp}
In the next two sections we show that there exists a family of $\l$ for which $i_*(\pi_1(\Q_\l))=ker(AJ_*)$. We do so by explicitly constructing loops of quadratic differentials corresponding to the generators of $ker(AJ_*)$; in this section we construct the necessary (square) transpositions, and in the next section the null $\rho_r$ and $i$-commutators. For simplicity of notation we refer to the elements we construct as being in $\pi_1(\Q_\l)$; however, they will be unchanged under $i_*$.  
 
We say that $\pi_1(\Q_\l^0, (M,q))$ \textit{contains all transpositions between} $p_i, p_j$ if for all edges $e \in E_{M,q}$ between $p_i$ and $p_j$, $\sigma_e$ or $\kappa_e \in i_*(\pi_1(\Q_\l^0, (M,q)))$ (depending on whether $p_i$ and $p_j$ are of the same weight or not). We say that $\pi_1(\Q_\l^0, (M,q))$ \textit{contains all transpositions} if the above is true for all pairs of zeroes of $q$. In Lemmas \ref{lem:+-} - \ref{lem:sym} we construct a variety of transpositions in hyperelliptic strata, and in the final proposition of the section we show that for some $\l$ the transpositions constructed in the four lemmas suffice to generate all transpositions in $\pi_1(\Q_\l^0, (M,q))$. 

For any hyperelliptic quadratic differential, $(M,q)$, we let $\tau$ denote the hyperelliptic involution on $M$ and $\pi$ a projection of $M$ to $\P^1$. We will use $\pi(M,q)$ to denote the the projection of $q$ to $\P^1$. Notice that any zero of a hyperelliptic quadratic differential at a branch point of $\pi$ must be of even order.  Unless specified otherwise we do not assume that the $k_1,...,k_n$ in $\Q_g(k_1^{m_1},k_2^{m_2},...,k_n^{m_n})$ are distinct. 

\begin{lemma} \label{lem:+-}
Let $(M_0,q_0) \in \Q_\l^0=\Q_g^0(k_1^2,...,k_m^2,k_{m+1},...,k_n)$  be a hyperelliptic quadratic differential, and label the zeroes of $q_0$ by $p_1^\pm,....,p_m^\pm, p_{m+1},...,p_n$. Suppose exactly $p_{m+1},...,p_n$ are at branch points of $M_0$ and at least one of $k_1,...,k_m$, say $k_l$, is even. Let $e$ be an edge between $p_i^+, p_j^+$, $1 \le i < j \le m$, such that $e \cap \tau(e) = \emptyset$. Then $\pi_1(\Q_\l^0, (M_0,q_0))$ contains an element corresponding to the (square) transposition of $p_i^+, p_j^+$ along $e$. Analogous statements are true for $p_i^+, p_j^-$, and $p_i^-, p_j^-$. 
\end{lemma}

\begin{figure}
\begin{center}
\includegraphics[width=.4 \textwidth]{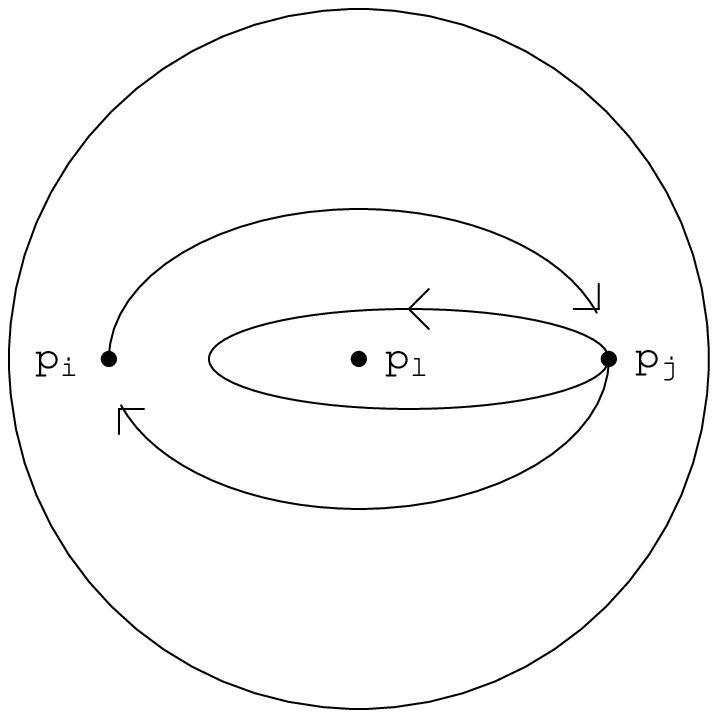}
\ \ \ \ 
\includegraphics[width=.4 \textwidth]{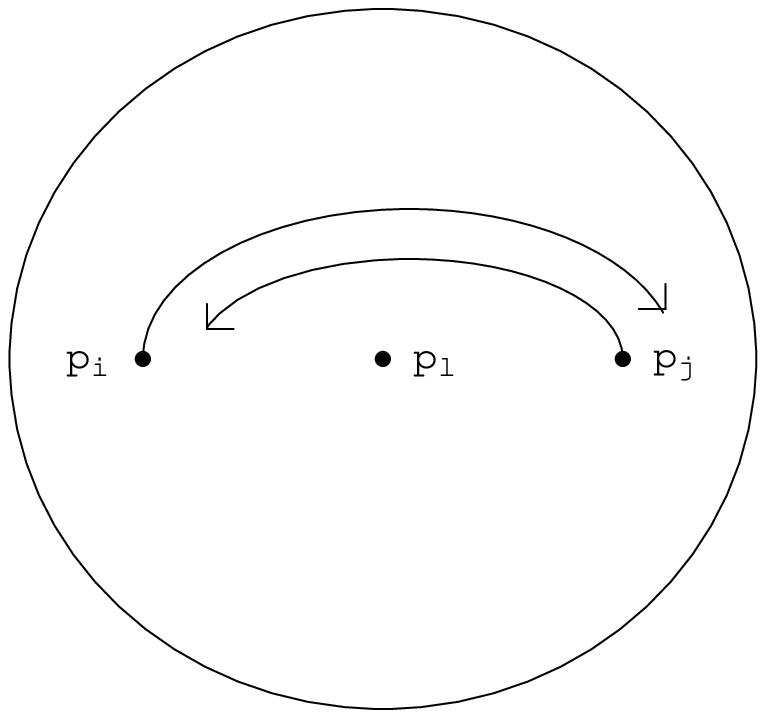}
\end{center}
\caption{Homotopic paths along which to transpose $p_i, p_j$, of odd order when $p_l$ is of even order.}
\label{fig:comp}
\end{figure}

\begin{proof}[Sketch of proof]
To prove the lemma we construct paths in the hyperelliptic locus of $\Q_\l$, from $(M_0,q_0)$ to quadratic differentials with discs that satisfy the conditions of Lemma \ref{lem:col} (as in Figure \ref{fig:coll}), and we then apply Proposition \ref{prp:tra}. In the case where at least one of $k_i, k_j$ is even, this proves the lemma. In the case where both $k_i$ and $k_j$ are odd the path we create moves $p_i, p_j$ near an even zero, $p_l$, and we create a (square) transposition of $p_i, p_j$ mod $p_l$. However, since $k_l$ is even we also have a square transposition of $p_j$ and $p_l$, and the composition of the two gives us the desired (square) transposition of $p_i$ and $p_j$ (see Figure \ref{fig:comp}). 
\end{proof}
\begin{proof}
Let $e: [0,1] \into M_0$, $e(0)=p_i^+, e(1)=p_j^+$, and let $\tilde{e}=\pi \cdot e:[0,1] \to \P^1$ be the projection of $e$ to $\P^1$. By abuse of notation we also use $e, \tilde{e}$ to mean the images of $[0,1]$ under the maps $e$ and  $\tilde{e}$.  

We first assume at least one of $k_i, k_j$ is even. Define: $$(\P^1, \tq_0)=\pi(M_0,q_0) \in \Q_{\tilde{\l}} = \Q_0(k_1,...,k_m,\frac{k_{m+1}-2}{2},...,\frac{k_n-2}{2}, -1^{2g+2-(n-m)})$$ and label the singularities of $\tq_0$ by $\tp_1, \tp_2,...,\tp_{2g+2+m}$, such that $\pi(p_s^\pm)= \tp_s$ for $1 \le s \le m$ and  $\pi(p_t)= \tp_t$ for $m < t \le n$. Since any branch zero is of even order, $(n-m) \le 2g-2$ and the number of single poles of $\tq$ is always positive. By Proposition \ref{prp:sphere} we may construct a quadratic differential on $\P^1$ with arbitrary zeroes; thus, define $(\P^1, \tq_t) \in \Q_{\tilde{\l}}$ to be the quadratic differential with zeroes at $\tp_1,...,\tp_{i-1},\tilde{e}(t), \tp_{i+1},...,\tp_{2g+2+m}$. Let $(M_t, q_t) \in \Q_\l^0$ be the double cover of $(\P^1, \tq_t)$, ramified at $\tp_{m+1}, ... ,\tp_{2g+2+m}$.  

As $t$ approaches 1, $e(t)$ approaches $p_j^+$, and there exists $t_0$ close to 1 such that $(M_{t_0}, q_{t_0})$ satisfies the conditions of Lemma \ref{lem:col}. Then we may apply Proposition \ref{prp:tra} to obtain a (square) transposition, $T$, of $p_j, e(t_0)$, in $\pi_1(Q_\l^0, (M_{t_0}, q_{t_0}))$, along a trajectory homotopic to $e([t_0, 1])$. Let $E$ be the image of $[0, t_0]$ in $\Q_\l^0$ under $t \mapsto (M_t, q_t)$. Along $E^{-1} \circ T \circ E$ $p_i^+$ and $p_j^+$ (square) transpose along $e$. Further, $p_i^-$ follows a trivial path (since $e \cap \tau(e) = \emptyset$) and all other zeroes stay constant, so this gives us the desired (square) transposition. 

Now suppose both $k_i$ and $k_j$ are odd.  Let $e$ be as above and recall that $k_l$ is even. let $\tilde{e}':[0,1] \into \P^1$ be an edge from $\tp_l$ to $\tp_j$ such that $\tilde{e}' \cap \tilde{e} = \tp_j$, and let $e':[0,1] \into M$ be the edge between $p_j^+$ and either $p_l^+$ or $p_l^-$ that projects to $\tilde{e}'$. Since $\tilde{e}'$ does not intersect itself $e' \cap \tau(e') = \emptyset$, and also by construction $e' \cap e = p_j^+$.

Define $(\P^1, \tq_t') \in \Q_{\tilde{\l}}$ to be the quadratic differential with zeroes at $\tp_1,...,\tp_{i-1},$ $\tilde{e}(t), \tp_{i+1},...,\tp_{l-1},\tilde{e}_t', \tp_{l+1},...,\tp_{2g+2+m}$, and let $(M_t', q_t') \in \Q_\l^0$ be the double cover of $(\P^1, \tq_t')$. Again there will be some $t_0$ such that $(M_{t_0}', q_{t_0}')$ satisfies the conditions of Lemma \ref{lem:col}. (Some slight deformation of $e$ and $e'$ may be required to get $e(t), e'(t)$ and  $p_j^+$ co-linear). Then by Proposition \ref{prp:tra}, $\pi_1(\Q_\l, (M_{t_0}', q_{t_0}'))$ contains a transposition of $e(t_0)$ and $p_j$ mod $e'(t_0)$ along $e([t_0, 1])$. However, since $e'(t_0)$ is a zero of even order there exists a square transposition of $e'(t_0)$ and  $p_j$ along $e'([t_0, 1])$, also in $\pi_1(\Q_\l, (M_{t_0}, q_{t_0}))$.  The composition of this square transposition and the transposition of $e(t_0)$ and  $p_j$ mod $e'(t_0)$ is shown in Figure \ref{fig:comp}, and gives a true transposition of $e(t_0)$ and $p_j$ along $e([t_0, 1])$. The rest of the proof follows as in the even case.  
\end{proof}

Notice that any geodesic, $e$, between two branch points of a hyperelliptic quadratic differential will have a `twin' geodesic, $\tau(e)$, parallel and of the same length. We would like to prove a lemma similar to $\ref{lem:+-}$ for branch zeroes; however, colliding two branch zeroes along $e$ will also cause the length of $\tau(e)$ to go to zero and a homology cycle to collapse. We deal with this by leaving the hyperelliptic locus of $\Q_\l$ shortly before colliding the two zeroes, and showing that for a generic element outside the hyperelliptic locus,  the trajectories corresponding to $e$ and $\tau(e)$ will be of different lengths.  

In particular, for any $(M,q)$ pick $x \in M$ and notice that for any loop, $\gamma$ on $M \bs P$ based at $x$ and $v \in T_xM$, parallel transport along $\gamma$ takes a vector $v$ either to $v$ or $-v$. This gives us a homomorphism: $$hol: \pi_1(M \bs P,x ) \to \Z_2$$ Following Masur and Zorich in \cite{MZ1}, this gives us a double cover of $M$, $\hM$, endowed with an induced flat metric that by construction has trivial holonomy and thus defines an Abelian differential, $\hq$ on $\hM$, with zeroes $\hat{P}$. $\hM$ is the \textit{canonical double cover} of $M$, and has a natural involution, $\phi$. Let $H_1^-(\hM, \hat{P};\Z)$ be the subspace of $H_1(\hM, \hat{P};\Z)$ that is anti-invariant with respect to $\phi$. Any $e \in H_1(M, P; \Z)$ has a double cover $\he' \cup \he''$, and we define $\he =\he'-\he'' \in H_1^-(\hM, \hat{P};\Z)$. We say that $e_1, e_2 \in H_1(M, P; \Z)$ are $\hh omologous$ if $\he_1$ is homologous to $\he_2$.

We may define a basis for $H_1^-(\hM, \hat{P}; \Z)$ by lifting a basis of $H_1(M, P; \Z)$, $\{e_i\}$, to the corresponding $\{\he_i\}$. Define the \textit{period} of $\he$ to be $\int_{\he} \hq$. Notice that locally elements of $\Q_\l$ will have double covers in the same stratum of Abelian differentials, and thus we may identify copies of $H_1^-(\hM, \hat{P}; \Z)$ over nearby surfaces. It is well known that the periods of a basis for $H_1^-(\hM, \hat{P}; \Z)$ give local coordinates on $\Q_\l$. Thus if $e_1$ and $e_2$ are not $\hh$omologous, the homology classes of $\he_1$, $\he_2$ are independent in $H_1(\hM, \hat{P}, \Z)$, and $| \int_{\he_1} \hq|=| \int_{\he_2} \hq|$ only on a set of measure zero.
Since $|e| = |\int_e q | = |\int_{\he'} \hq |=\frac12 |\int_{\he} \hq |$, the same will be true of $e_1$ and $e_2$ -  if they are non-$\hh$omologous they are of the same length only on a set of measure zero. 

\begin{lemma} \label{lem:cdc}
Let $(M_0,q_0)$ and $\l$ be as in Lemma \ref{lem:+-}, with some $k_i$ odd, at least two branch zeroes, and $g > 0$. Let $e^+$ be a trajectory between two of the branch zeroes of $q_0$ such that the trajectory does not contain any other zeroes of $q_0$, and let $e^-:=\tau(e^+)$. Then there exists some open set, $U \subset Q_\l$, containing $(M_0,q_0)$, such that for almost every element of $U$ these trajectories are not of the same length. 
\end{lemma}

\begin{proof} 
 As noted above, if $e^+$ and $e^-$ are not $\hh$omologous in $H_1^-(\hM_0, \hat{P}_0; \Z)$ then they are generically not of the same length. Thus, to prove the lemma we suppose that $\he := \he^+ - \he^-$ is homologous to zero and show a contradiction. 

If $\he$ is homologous to zero then it is a separating curve. Let $\tM_0$ be the surface with boundary obtained by cutting $M_0$ along $e^+ \cup e^-$, which together form a non-trivial cycle. Since $g > 0$, $\tM_0$ is connected.  If any element of $\pi_1(\tM_0, P_0)$ mapped to $-1$ under $hol:\pi_1(\tM_0, P_0) \to \Z_2$ then $\hM_0$ cut at $\he$ would still be connected, which implies that every element of $\pi_1(\tM_0,P_0)$  maps to 1 under $hol$. But by assumption one of the zeroes of $P_0$ is odd and a small loop around it will map to $-1$ under $hol$, giving us the desired contradiction. 
\end{proof}

Now we can prove the analog of Lemma \ref{lem:+-} for branch zeroes. 

\begin{lemma} \label{lem:br}
Let $(M_0,q_0)$, $\l$ and $\tilde{\l}$ be as in Lemma \ref{lem:+-}, with at least one $k_i$ odd, and let $e$ be any edge between two branch zeroes of $q_0$, $p_i$ and $p_j$, $m < i<j \le n$. Then $\pi_1(\Q_\l^0, (M_0,q_0))$ contains a (square) transposition of $p_i, p_j$ along $e$. 
\end{lemma}

\begin{proof}
Let $e:[0,1] \into M_0$, $e(0)=p_i$, $e(1)= p_j$, be the edge between $p_i, p_j$, and let $\tilde{e}=\pi \circ e$ be its projection to $\P^1$. Let $(\P^1, \tq_0)=\pi(M_0,q_0) \in \Q_{\tilde{\l}}$ and label the zeroes of $(\tM_0, \tq_0)$ as $\tp_1, \tp_2,..., \tp_{2g+2+n}$. Define $(\P^1, \tq_t)$ to be the element of $\Q_{\tilde{\l}}$ with zeroes at $\tp_1, \tp_2,...,\tp_{i-1}, \tilde{e}(t), \tp_{i+1},...,\tp_{2g+2+n}$. Taking a double cover of $(\tM_t, \tq_t)$ branched at its last $2g+2$ zeroes we get an element of $\Q_\l$, which we call $(M_t,q_t)$. 

As $t$ approaches 1, $e(t)$ approaches $p_j$ and all other $p_k$ are fixed. As in Lemma \ref{lem:+-} we wish to show there exists a subset of some $M_t$ containing $e(t)$ and $p_j$, and satisfying the conditions of Lemma \ref{lem:col}. In fact, because $e(t)$ and $p_j$ are both of even order we may actually relax the conditions of Lemma \ref{lem:col}: if $e(t)$ and $p_j$ are $\delta$ apart it suffices to show that all other saddle connections on $M_t$ are of length greater than $\delta$. (We see this by proving Lemma \ref{lem:col} in the same way but only for the case of two even zeroes, as in the second pair of discs in Figure \ref{fig:coll}). This condition is not initially satisfied because on $(M_t,q_t)$ there are two short trajectories of the same length running between $e(t)$ and $p_j$. However, Lemma \ref{lem:cdc} implies that there exists $(M,q)$ arbitrarily close to any $(M_t,q_t)$ such that the trajectories in the homotopy classes of $e([t,1])$ and $\tau(e([t,1]))$ are of different lengths. Because of the involution on $(M_t, q_t)$ we may assume the trajectory on $(M,q)$ in the homotopy class of $e([t,1])$ is shorter and thus $(M,q)$ satisfies the (modified) conditions of Lemma \ref{lem:col}. Then there exists a (square) transposition of $e(t), p_j$ along $e([t,1])$, based at $(M,q)$. The rest of the argument follows as in Lemma \ref{lem:+-}.
\end{proof}

\begin{lemma}  \label{lem:sym}
Let $(M_0, q_0)$ and $\l$ be as in Lemma \ref{lem:+-}, and let $b$ be a branch point of $M_0$ that is also a regular point of $q_0$. Let $e_b$ be an edge between $b$ and $p_i^+$, $1 \le i\le m$, such that $\tau(e_b) \cap e_b = b$. Then $e := e_b \cup \tau(e_b)$ is an edge between $p_i^+$ and $p_i^-$, and  $\pi_1(\Q_\l^0, (M_0,q_0))$ contains an element corresponding to the transposition of $p_i^+, p_i^-$ along $e$. 
\end{lemma}

\begin{proof}
Let $e'$ be a small deformation of $e$ to one side of $b$ that is still an edge. We can always choose $e'$ so that $e' \sim e \sim \tau(e')$ and $e' \cap \tau{e'}=p_i^\pm$. As in the previous two lemmas, let $(\P^1, \tq_0)$ be a projection of $(M_0,q_0)$ to $\P^1$, with singularities $\tp_1,...,\tp_{2g+2+m}$, and let $\tilde{e}'=\pi(e')$. Define $(\P^1, \tq_t)$ to be the quadratic differential with zeroes at $\tp_1,..., \tp_{i-1}, \tilde{e}'(t), \tp_{i+1},..., \tp_{2g+2+m}$, and let $(M_t,q_t)$ be its double cover. $(M_0,q_0)=(M_1,q_1)$ and $T:[0,1] \to \Q_\l$, $t \mapsto (M_t, q_t)$ is the desired element of $\pi_1(\Q_\l, (M_0,q_0))$. 
\end{proof}

The previous four lemmas will allow us to show that certain hyperelliptic strata contain all of their transpositions, but first we need one more technical lemma. Recall that the maximal number of faces a planar graph with $n$ vertices can have is $2n-4$. (Such a graph will be a triangulation, with $3n-6$ edges.) By removing edges from such a graph it is always possible to construct a planar graph with $n$ vertices and fewer than $2n-4$ faces. (As earlier in the paper, we do not allow graphs to have double edges or loops.)

\begin{lemma} \label{lem:tgr}
Let $\Gamma$ be a planar graph with $n$ vertices and $f$ faces, $f \le 2n-4$. Then we may associate to each face of $\Gamma$ a pair of vertices adjacent to the face such that the same pair is not associated to more than 1 face. 
\end{lemma}

\begin{proof}
It suffices to associate a unique adjacent edge to each face of $\Gamma$, since this is equivalent to associating the two vertices adjacent to the edge to the face. Pick any face, $F_1$, of $\Gamma$ and any edge, $e_1$, adjacent to $F_1$. $e_1$ is adjacent to one other face, which we call $F_2$. $F_2$ has at least two other possible adjacent edges, so again pick any  edge $e_2 \ne e_1$ to associate to $F_2$, and let $F_3$ the other face adjacent to $e_2$. Continue this process until one of two things happens. Either an edge is associated to each face, or at stage $k$ $F_k=F_i$ for $1 \le i \le k$ and $k \ne f$. In the first case we are done, and in the second the remaining faces of the graph do not have any adjacent edges that have been associated to any other face, so we pick any face, call it $F_{k+1}$ and resume the process.   
\end{proof}

\begin{prop} \label{prp:hy2}
Suppose that $\Q_\l=\Q_g(1^{2n}, k_1,...,k_m)$ is such that $2n \ge g+5$, $k_1,..,k_m$ are all even, and there exists $i,j$ such that $k_i=k_j$. Then for any $(M,q) \in \Q_\l$, $\pi_1(\Q_\l, (M,q))$ contains all transpositions.  
\end{prop}

\begin{figure}
\begin{center}
\includegraphics[width=.4 \textwidth]{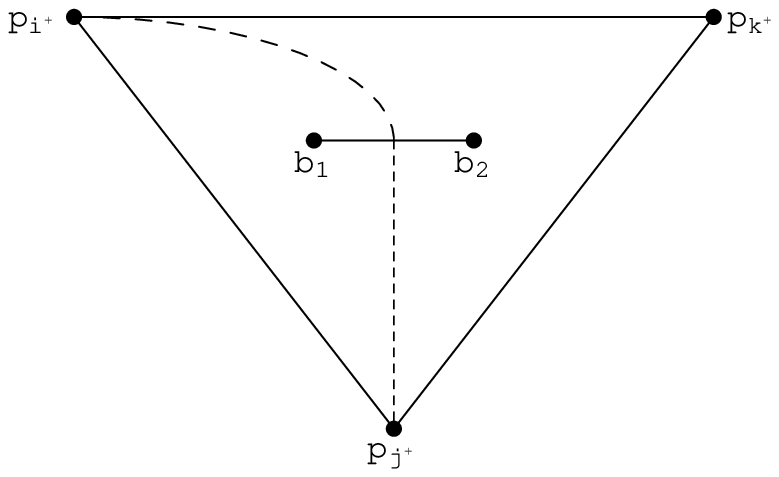}
\ \ \ \ 
\includegraphics[width=.4 \textwidth]{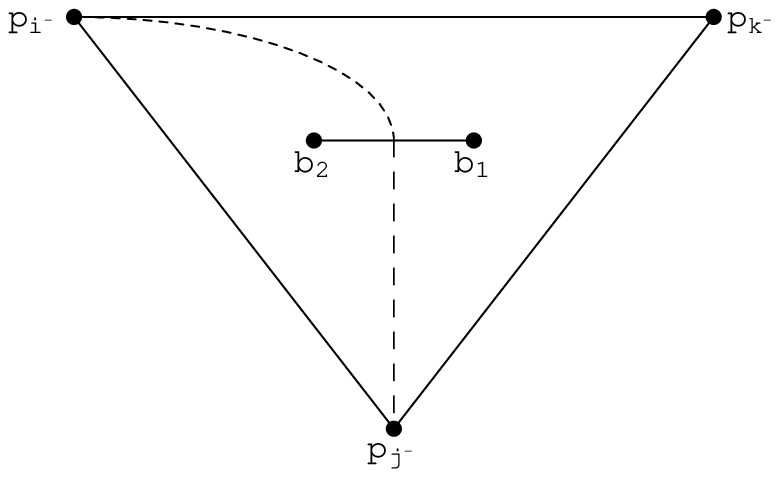}
\end{center}
\caption{The edges $e_k$ and $\tau(e_k)$ between $p_i^+, p_j^-$ and $p_i^-, p_j^+$. $b_1$ and $b_2$ denote branch points, and the line between them denotes a branch cut.}
\label{fig:branch}
\end{figure}

\begin{proof}[Sketch of Proof]
Our plan is to construct a hyperelliptic $(M,q)$ with higher order zeroes at branch points, and then put a graph on $M$ with the single zeroes of $q$ at the vertices and at most one higher order zero in each face. We use Lemmas \ref{lem:+-} and \ref{lem:sym} to construct the transpositions associated to each edge of the graph, which via Theorem \ref{thm:copeland} gives us all (square) transpositions of single zeroes with each other and with higher order zeroes. We then apply Lemma \ref{lem:br} to get (square) transpositions of the higher order zeroes with each other. 
\end{proof}
\begin{proof}
By Theorem \ref{thm:c1} $\Q_\l$ is connected so we need not consider connected components. Without loss of generality suppose $k_1=k_2$. Put a graph $\tilde{\Gamma}$ on $\P^1$ with $n$ vertices and $g+1$ faces. Since $n \ge \frac{g+5}{2}$, $g+1 \le 2n-4$, and it is always possible to construct such a graph. Mark 2 points in each of $g$ faces of $\tilde{\Gamma}$, and 3 points in the $(g+1)$st. By Proposition \ref{prp:sphere} there exists $\tq$ on $\P^1$ such that the vertices of $\tilde{\Gamma}$ are zeroes of $\tq$ of order $1$, $m-1$ marked points are zeroes of $\tq$ of order $k_1, \frac{k_3-2}{2}, \frac{k_4-2}{2},...,\frac{k_m-2}{2}$, and the $2g+3-m$ remaining marked points are poles of $\tq$ of order 1. Since $m < 2g+2$, $\tq$ has at least 2 poles and we assume the $g+1st$ face contains the zero of order $k_1$ and 2 poles.  Denote the $n$ single zeroes of $\tq$ by $p_1,...,p_n$. To each face of $\tilde{\Gamma}$ we associate a pair of vertices as in Lemma \ref{lem:tgr}. 

Take a double cover of $(\P^1, \tq)$ ramified at the $2g+2$ marked points in the faces of $\tilde{\Gamma}$, minus the zero of order $k_1$. Assume each branch cut is between two points in the same face and is contained in that face. This gives us $(M,q) \in \Q_g(1^{2n}, k_1,k_2,....,k_m)$. $(M,q)$ has two copies of $\tilde{\Gamma}$ embedded into it, $\tilde{\Gamma}^+$ and $\tilde{\Gamma}^-$; denote their vertices by $p_1^+,...,p_n^+$ and $p_1^-,...,p_n^-$, with $\tau(p_l^+)=p_l^-$, $1 \le l \le n$. By associating a pair of vertices to each face of $\tilde{\Gamma}$ we have associated a quadruplet of vertices, $p_i^\pm$ and $p_j^\pm$, to each of the $g+1$ branch cuts of $M$. Construct a pair of edges, $e_k$ and $\tau(e_k)$, between $p_i^+, p_j^-$ and $p_j^+, p_j^-$ through the $k$th branch cut, $1 \le k \le g$, as in Figure \ref{fig:branch}. Notice that we may always construct $e_k$ and $\tau(e_k)$ so that they do not intersect. For the $(g+1)$st branch cut, both branch points are regular and we construct edges, $e_{g+1}$ and $e_{g+1}'$ between $p_i^\pm$ and $p_j^\pm$, as in Lemma \ref{lem:sym}. 
$\tilde{\Gamma}^+ \cup \tilde{\Gamma}^- \cup e_1 \cup \tau(e_1) \cup ... \cup e_g \cup \tau(e_g)\cup e_{g+1} \cup e_{g+1}'$ gives us a connected graph $\Gamma$ on $M$ with $2g+2$ faces and two faces associated to each branch cut. This graph will be a 2-cell embedding and by construction satisfies all of the hypotheses of Theorem \ref{thm:copeland}. Exactly one branch point is in each of faces associated to the first $g$ branch cuts. For the two faces, $F$ and $F'$, associated to the $g+1$st branch cut, $F'=\tau(F)$ so again there is one zero of order $k_1(=k_2)$ in each of $F, F'$.

By Lemma \ref{lem:+-} we may construct an element of $\pi_1(\Q_\l, (M,q))$ corresponding to the transposition of $p_i^+, p_j^-$ along both $e_k$ and $\tau(e_k)$ for $1 \le k \le g$. Further, each edge $e^+$ of $\tilde{\Gamma}^+$ is contained entirely on one sheet of $M$, so $e^+ \cap \tau(e^+) = \emptyset$ and again by Lemma \ref{lem:+-} $ \sigma_{e^+}, \sigma_{e^-} \in \pi_1(\Q_\l, (M,q))$.   
Lemma \ref{lem:sym} gives us transpositions associated to $e_{g+1}$ and $e_{g+1}'$. Consequently for every edge $e$ of $\Gamma$,  $\sigma_e \in \pi_1(\Q_\l, (M,q))$. Since there is at most one higher order zero in each face, Theorem \ref{thm:copeland} implies $\pi_1(\Q_\l, (M,q))$ contains all transposition of the $2n$ zeroes of order 1 with each other, and all square transpositions of zeroes of order 1 with zeroes of order $k_1,...,k_m$.  

Now pick any $(\P^1, \tq') \in \Q_0(1^n,  \frac{k_1-2}{2},...,\frac{k_m-2}{2}, -1^{2g+2-m})$, and let $(M',q') \in \Q_\l$ be its double cover, ramified at the zeroes of $\tq'$ not of order 1. Then all of the even zeroes of $q'$ are at branch points of $M'$ and by Lemma \ref{lem:br} $\pi_1(\Q_\l, (M',q'))$ contains all (square) transpositions of zeroes of order $k_1,...,k_m$ with each other.  But $\pi_1(\Q_\l, (M',q'))$ is isomorphic to $\pi_1(\Q_\l, (M,q))$ so the same is then true for $\pi_1(\Q_\l, (M,q))$. This proves the proposition. 
\end{proof}

\section{Constructing the remaining generators}\label{sec:null}

We have shown that in some cases all transpositions and square transpositions of zeroes of $q$ are contained in $\pi_1(Q_\l^0, (M,q))$. However, for $\l$ where many but not all points are not of equal weight, Theorem \ref{thm:gen} and Corollary \ref{cor:gen2} imply that $ker(\pi_1(Sym_g^\l) \to H_1(M, \Z))$ is generated by transpositions, null $\rho_r$, and $i$-commutators. In some cases we can again use techniques of colliding and breaking apart points from Section \ref{sec:adj} to show that these last two types of elements are in $\pi_1(\Q_\l^0, (M,q))$.

Let $l_1,...,l_{2g}$ be standard generators of $\pi_1(M)$. Recall that in Sections \ref{sec:pi1} and \ref{sec:ker} we did not explicitly define the $\rho_{ir}$, we only stated that one of the generators of $B_\l$ must be a loop corresponding to $p_i$ moving around $l_r$. There are infinitely many choices of such a loop, differing by various products of transpositions. 

\begin{lemma} \label{lem:gross}
Let $\Q_{\l_1}=\Q_g(k_1,k_2,...,k_n)$ contain all of its transpositions, and suppose $\sum_{i=1}^l k_{i} = \sum_{j=l+1}^m k_{j}$. Further suppose $\Q_{\l_1}^0$ is adjacent to $\Q_{\l_2}^0=\Q_g^0(k_{1}+...+k_{l}, k_{l+1}+...+k_{m}, k_{m+1}, ...,k_n) \ni (M,q)$ and $\pi_1(\Q_{\l_2}^0, (M,q))$ contains all transpositions of the two newly formed points. Define $\a(1), \a(2),...,\a(l)$ to all equal 1, and $\a(l+1),...,\a(m)$ to all equal -1. Then for any $(M',q') \in \Q_{\l_1}$, for any choice of $\rho_{1r},...,\rho_{nr}$, for any $i_1,...,i_m$ such that $k_{i_s}=k_s$, $1 \le s \le m$, and for any $\sigma \in S_m$, $\rho_{i_{\sigma(1)}r}^{\a(\sigma(1))}\cdot...\cdot\rho_{i_{\sigma(m)}r}^{\a(\sigma(m))} \in \pi_1(\Q_{\l_1}^0, (M',q'))$ .  
\end{lemma}

\begin{figure}
\begin{center}
\includegraphics[width=.3 \textwidth]{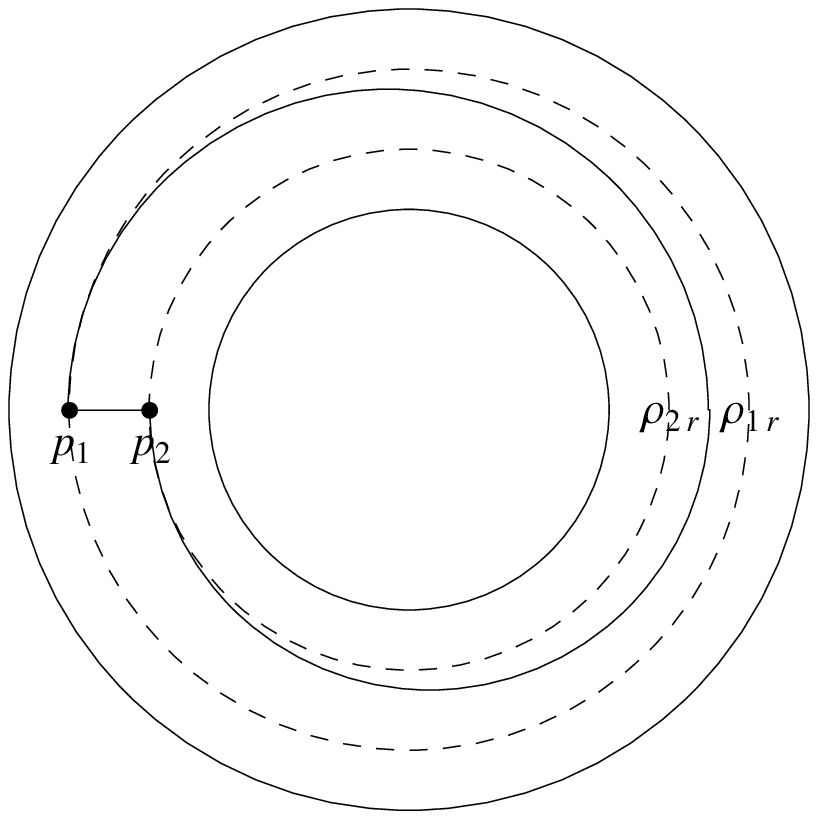}
\ \ \ \
\includegraphics[width=.3 \textwidth]{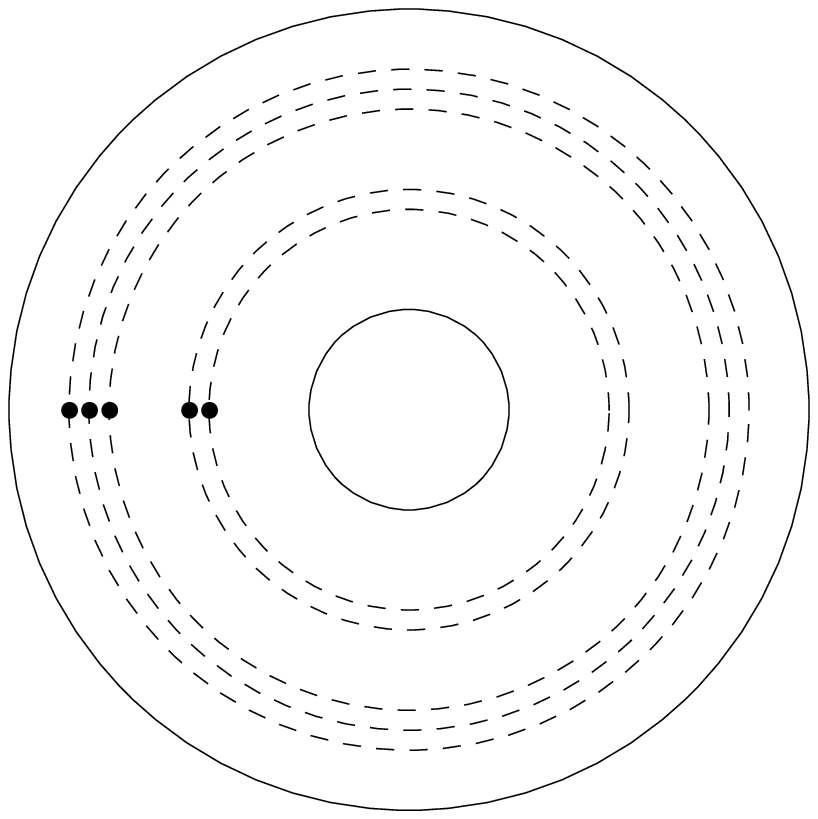}
\end{center}
\caption{On the left, transposing $p_1, p_2$ along both the solid curves is equivalent to $\rho_{1r}^\pm \rho_{2r}^\mp$. On the right, surgering give us a number of $\rho_{ir}$ from the original two.}
\label{fig:null}
\end{figure}

\begin{proof}
To prove the lemma we first explicitly define $\rho_{1r},..., \rho_{nr}$ and show that under our choice of definition, $\rho:=\rho_{1r}...\rho_{lr}\rho_{(l+1)r}^{-1}...\rho_{mr}^{-1} \in \pi_1(\Q_\l, (M',q'))$. We then pick any other choice of generators, $\rho'_{1r},..., \rho'_{nr}$, any $i_1,...,i_m$ (not necessarily distinct) such that $k_{i_s}=k_s$, and any permutation $\sigma \in S_m$. We define $\rho':=\rho_{i_{\sigma(1)}r}^{'\a(\sigma(1))}\cdot...\cdot\rho_{i_{\sigma(m)}r}^{'\a(\sigma(m))}$. We show that $\rho'$ differs from $\rho$ only by a product of transpositions and is thus also in $\pi_1(\Q_\l, (M',q'))$.

Label the two zeroes of $q$ of order $k_1+k_2+...+k_l$ as $p_1$ and $p_2$. For any $r$ pick two edges, $e$ and $e'$, between $p_1, p_2$ such that $e \cup e'$ is homotopic to $l_r$ and $e \cap e'=\{p_1, p_2\}$, as in the left of Figure \ref{fig:null}. By assumption $ \sigma_{e}, \sigma_{e'} \in \pi_1(\Q_{\l_2}^0, (M,q))$. Further, defining $\rho_{1r}$ and  $\rho_{2r}$ as in the left of Figure \ref{fig:null}, $\sigma_{e} \sigma_{e'}$ is homotopic to $\rho_{1r}\rho_{2r}^{-1}$. Thus $\rho_{1r}\rho_{2r}^{-1} \in \pi_1(\Q_{\l_2}^0, (M,q))$.

 Now pick an explicit loop of surfaces in the homotopy class $\rho_{1r} \rho_{2r}^{-1}$, $\eta: [0,1] \to \Q_{\l_2}^0$, such that $\eta(0)=\eta(1)=(M,q)$. Surger every $\eta(t)$ by some continuously varying $\theta_t, \delta_t$, possibly multiple times, to get a loop of surfaces $\eta':[0,1] \to \Q_{\l_1}$, and define $(M',q') = \eta'(0)$. Let $p_{1}(t),...,p_{m}(t)$ be the zeroes of $\eta'(t)$ of order $k_{1},...,k_{m}$ formed from the surgery. We may then define $\rho_{1r}$ to be the path followed by $p_{1}(t)$, and so forth (as in the right of Figure 8). Under these definitions of the $\rho_{sr}$, $1 \le s \le m$, we have created  $\rho:=\rho_{1r}...\rho_{lr}\rho_{(l+1)r}^{-1}...\rho_{mr}^{-1} \in \pi_1(\Q_{\l_1}, (M',q'))$. 

Now let $\sigma$ be an element of $S_m$ and consider $\rho_{\sigma}:=\rho_{\sigma(1)r}^{\a(\sigma(1))}...\rho_{\sigma(m)r}^{\a(\sigma(m))}$.  Since $\rho_\sigma \cdot \rho^{-1}$ is supported on a subset of $M'$ that is a disc, it may be written as a product of (square) transpositions. We have assumed $\pi_1(\Q_\l, (M',q'))$ contains all transpositions, so $\rho_\sigma \cdot \rho^{-1} \in \pi_1(\Q_\l, (M',q'))$ and therefore $\rho_{\sigma} \in \pi_1(\Q_\l, (M',q'))$. 

Similarly define $\rho'_{1r}, \rho'_{2r},...,\rho'_{nr}$ be a different choice of generators, and define $\rho':=\rho'_{1r}...\rho'_{lr}\rho^{'-1}_{(l+1)r}...\rho^{'-1}_{mr}$. Again $\rho' \cdot \rho^{-1}$ is supported in a disc on $M'$ and therefore $\rho' \in \pi_1(\Q_\l, (M',q'))$.

 Finally suppose we have $p_{i_1},....,p_{i_m}$ (not necessarily distinct) such that $k_{i_s}=k_s$, $1 \le s \le m$, and define $\rho_i:=\rho_{i_1}...\rho_{i_l}\rho_{i_{l+1}}^{-1}...\rho_{i_m}^{-1}$. Either $s=i_s$ and $\rho_{sr}^\pm \rho_{i_sr}^\mp=1$ or $\rho_{sr}^\pm \rho_{i_sr}^\mp$ may be written as a product of two transpositions, as in the first paragraph of this proof. 
Then the following is in $\pi_1(\Q_\l, (M',q'))$: 
$$\beta:=\rho_{1r}\rho_{i_1 r}^{-1}\rho_{2r}\rho_{i_2 r}^{-1}...\rho_{lr}\rho_{i_l r}^{-1}\rho_{(l+1)r}^{-1}\rho_{i_{l+1}r}...\rho_{m}^{-1}\rho_{i_{m}r}$$ 
As above $\rho_i \beta \rho^{-1}$ is supported in a disc and therefore in $\pi_1(\Q_\l, (M',q'))$. This then implies $\rho_i \in \pi_1(\Q_\l, (M',q'))$.

A combination of the above implies the lemma for the specific choice of $(M',q')$, surgered from $(M,q)$. However, since fundamental groups with different base points are isomorphic, the same will be true for any element of $\Q_{\l_1}$. 
\end{proof}

Lemma \ref{lem:gross} shows that under certain adjacency conditions it is possible to construct any null $\rho_r$ involving points of certain weights. When $l,m=1$ the lemma implies that if $\pi_1(\Q_{\l_1}^0, (M',q')$ contains all of its transpositions, it also contains all null $\rho_r$ involving two points of equal weights. 

\begin{prop} \label{prp:null}
 Let $g > 2$ and $\l_1=(1^a,k_1,...,k_n)$ with $a >$ max$\{g+4,k_1,..., k_n \}$, all $k_i$ even, and some $k_i=k_j$. Then for any $(M',q')$, all null $\rho_r$ are contained in $\pi_1(Q_{\l_1}, (M',q'))$.
\end{prop}

\begin{proof}
By Theorem \ref{thm:c1} there is only one connected component of $\Q_{\l_1}$. 

Let $(M',q') \in \Q_{\l_1}$ with $p_1,...,p_a$ the zeroes of $q$ of order 1. By Proposition \ref{prp:hy2} $\pi_1(\Q_{\l_1}, (M',q'))$ contains all of its transpositions; by Lemma \ref{lem:gross} it contains any $\rho_{jr} \rho_{lr}^{-1}$, $1 \le j,l \le a,1 \le r \le 2g$. 

For $k_i \ne 2$, let $\Q_{\l_2}=\Q_g(1^{a-k_i}, k_1,...,k_i^2,...,k_n)$ and notice that by Lemmas \ref{lem:two} and \ref{lem:three}, $\l_1 > \l_2$ (this is not true if $k_i=2$). Since $k_i$ is even there exists a component of $\Q_{\l_2}$, $\Q_{\l_2}^0$, that contains a hyperelliptic element, $(M,q)$, for which both zeroes of order $k_i$ are at branch points. By Lemma \ref{lem:br}, $\pi_1(\Q_{\l_2}^0, (M,q))$ contains all transpositions of the two zeroes of order $k_i$, and by Lemma \ref{lem:gross} it contains each of the $2g$ null $\rho_r$ involving only those two zeroes. Therefore Lemma \ref{lem:gross} implies any null $\rho_r$ consisting of a single zero of order $k_i$ moving around $l_r$ and $k_i$ zeroes of order 1 moving around $l_r^{-1}$ is in $\pi_1(\Q_{\l_1}, (M',q'))$.

If $k_i=2$ we let $\Q_{\l_2}=\Q_g(1^{a-6},4^2, k_1,...,\hat{k}_i,...,.k_n)$ and note that $\l_1 > \l_2$. By the same argument as above there exists $(M,q) \in \Q_{\l_2}$  such that $\pi_1(\Q_{\l_2}, (M,q))$ contains each of the $2g$ null $\rho_r$ involving only the two zeroes of order 4. Thus we get the null $\rho_r$ involving a zero of order $2$ and $2$ zeroes of order 1 moving one way around $l_r$, and 4 zeroes of order 1 moving the other. We compose and cancel with null $\rho_r$ involving points of order 1 to get a null $\rho_r$ with a zero of order 2 moving around $l_r$ and two zeroes of order 1 moving around $l_r^{-1}$. 

This gives us any null $\rho_r$ consisting of a single zero of higher order moving one way around $l_r$ and zeroes of order 1 moving the other. These and the $\rho_{jr} \rho_{lr}^{-1}$ where $k_j=k_l=1$ generate all null $\rho_r$.
\end{proof}

Notice that there are strata of the form specified in Proposition \ref{prp:hy2} for which there exist null $\rho_r$ to which Lemma \ref{lem:gross} does not apply. For example, in $\Q_{10}(1^{16},20)$ we can have a null $\rho_r$ consisting of the point of order 20 moving one way around $l_r$ and 20 points of order 1 moving the other, but since there are not 20 distinct points of order 1 we cannot collide them to use the technique of Lemma \ref{lem:gross}.

Finally we would like to consider when $i$-commutators are contained in $\pi_1(\Q_\l)$. 

\begin{prop} \label{prp:comm}
Let $\l=(1^{a}, k_1,...,k_n)$, where $a$ and the $k_i$ are as in Proposition \ref{prp:null}. Then for any $(M,q) \in \Q_\l$, any $i$-commutator is in $\pi_1(Q_\l, (M,q))$.
\end{prop}

\begin{figure}
\begin{center}
\includegraphics[width=.3 \textwidth]{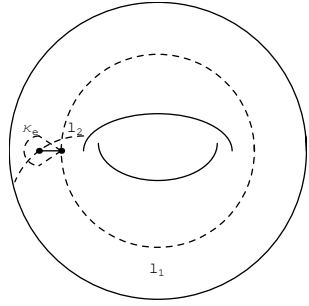}
\end{center}
\caption{$\rho_{11}, \rho_{22}$, and $\kappa_e$ defined on a torus with two marked points. The solid line is $e$.}
\label{fig:kappa}
\end{figure}

\begin{proof} 
Let $p_1, p_2,...,p_a$ be the zeroes of $q$ of weight 1, and $p_{a+1},...,p_{a+n}$ the zeroes of weight $k_1,...,k_n$. where $k_i$ is by assumption is less than $a$, $1 \le i \le n$. By Proposition \ref{prp:null} $\pi_1(\Q_\l, (M,q))$ contains all of its null $\rho_r$. Thus the following is also in $\pi_1(Q_\l, (M,q))$:
\begin{equation} \label{eq:rho1}
(\rho_{1r} ... \rho_{kr} \rho_{ir}^{-1})(\rho_{1s} ... \rho_{ks} \rho_{is}^{-1})(\rho_{1r}^{-1} ... \rho_{kr}^{-1} \rho_{ir}) (\rho_{1s}^{-1} ... \rho_{ks}^{-1} \rho_{is})
\end{equation}
For $i \ne j$ and $l_r \cap l_s = \emptyset$, $\rho_{ir}$ commutes with $\rho_{js}$. If $i \ne j$ but $l_r \cap l_s = 1$ then  $\rho_{ir} \rho_{js}$ =  $\rho_{js} \rho_{ir} \kappa_{e}$ where $\kappa_{e}$ is a square transposition of $p_i, p_j$, defined appropriately with respect to $\rho_{js}, \rho_{ir}$. For example, in Figure \ref{fig:kappa} we have two points going around $l_1, l_2$ on a torus and an explicitly defined $\kappa_e$, with $\rho_{11} \rho_{22} = \rho_{22} \rho_{11} \kappa_e$.

 Thus if $l_r \cap l_s = \emptyset$ (\ref{eq:rho1}) is equal to:
\begin{equation} \label{eq:rho2}
(\rho_{1r}...\rho_{kr}\rho_{1s} ... \rho_{ks}\rho_{1r}^{-1} ... \rho_{kr}^{-1}\rho_{1s}^{-1}...\rho_{ks}^{-1})(\rho_{ir}^{-1} \rho_{is}^{-1} \rho_{ir} \rho_{is})
\end{equation}
Notice that the maximal value for $n$ is $(4g-4-(g+5))/2$, and this implies that $g+5 \ge \frac{3 + \sqrt{9+8(2g+n-2)}}{2}$  for any $g$. The first of the two elements in parentheses in (\ref{eq:rho2}) is in the kernel of $AJ_*$, so we may apply Corollary \ref{cor:imv2} to show that it can be written as a product of transpositions. Proposition \ref{prp:hy2} says that $\pi_1(\Q_\l, (M,q))$ contains all of its transpositions; thus the element on the left is in $\pi_1(\Q_\l, (M,q))$. This in turn implies implies $\rho_{ir} \rho_{is} \rho_{ir}^{-1} \rho_{is}^{-1} \in \pi_1(\Q_\l, (M,q))$. 

If $l_r \cap l_s=1$ then (\ref{eq:rho1}) is equal to (\ref{eq:rho2}) except that the first element in parentheses will contain some additional $\kappa_e$'s. However it will still be in the kernel of $AJ_*$ and again $\rho_{ir} \rho_{is} \rho_{ir}^{-1} \rho_{is}^{-1} \in \pi_1(\Q_\l, (M,q))$. 

Similarly, for $p_l$ a zero of higher order and $e'$ an edge between $p_i, p_l$, we would like to show that $\kappa_{e'}^{-1} \rho_{ir}^{-1} \kappa_{e'} \rho_{ir}$ is in $\pi_1(\Q_\l, (M,q))$, so we consider the following: 

$$\kappa_{e'}^{-1}(\rho_{1r} ... \rho_{kr} \rho_{ir}^{-1})\kappa_{e'}(\rho_{1r}^{-1} ... \rho_{kr}^{-1} \rho_{ir}).$$
By assumption both $\kappa_{e'}^{\pm 1}$ and the elements in parentheses are in  $\pi_1(\Q_\l, (M,q))$, so the whole element is. Since $\kappa_{e'}$ commutes with $\rho_{1r},...,\rho_{kr}$, commuting the $\rho_{1r},...,\rho_{1k}$ as in 
(\ref{eq:rho2}) gives us the desired result. 

Finally, any product of square transpositions of two higher order zeroes is contained in $\pi_1(\Q_\l, (M,q))$ by Proposition \ref{prp:hy2}, because it is a product of transpositions. This covers all possible generators of $[\pi_1(M_{(n-i)}),\pi_1(M_{(n-i)})]$.
\end{proof}

\section{Conclusion} \label{sec:sum}
We summarize by answering the question of when $ker(AJ_*:\pi_1(Sym_g^\l) \to H_1(M\,Z))$ is equal to $im(i_*:\pi_1(\Q_\l) \to \pi_1(Sym_g^\l))$. 

\begin{theorem}Let $\l=(1^a,k_1,...,k_n)$ with $a >$ max$\{g+5,k_1,..., k_n \}$, all $k_i$ even, and some $k_i=k_j$. Then $im(i_*)=ker(AJ_*)$.
\end{theorem}

\begin{proof}
Proposition \ref{prp:triv} implies $im(i_*:\pi_1(\Q_\l) \to \pi_1(Sym_g^\l)) \subset ker(AJ_*)$. By Theorem \ref{thm:gen} and Corollary \ref{cor:gen2} $ker(AJ_*)$ is generated by transpositions, square transpositions, null $\rho_r$ and in some cases $i$-commutators. Proposition \ref{prp:hy2}, Proposition \ref{prp:null}, and Proposition \ref{prp:comm} show that all of these elements are in $im(i_*:\pi_1(\Q_\l) \to \pi_1(Sym_g^\l))$. \end{proof}

In \cite{C}, Copeland shows a similar result for $g > 2$ and $\l=(1^{4g-4})$. His techniques are somewhat different and rely on the fact that in the top stratum one may interpolate two quadratic differentials and expect the result to be in the same stratum. 

Thus, for certain $\l$ we have constructed $i_*(\pi_1(\Q_\l))$. Of course, we are actually interested in $\pi_1(\Q_\l)$ and would thus like to determine the kernel of $i_*$.  However, it may be difficult to say anything about this kernel. 

For example, let $\l=(1^{4g-4})$. $Sym_g^\l$ and $Pic_g^{4g-4}$ are both smooth connected complex varieties and $AJ:Sym_g^\l \to Pic_g^{4g-4}$ is a dominant morphism of varieties. A fiber of $AJ$ is a copy of $\P^{3g-4}$ with a codimension 1 subset removed, so it is connected and has at least one smooth point. 

Then there is a non-empty Zariski open set $U \subset Pic_g^{4g-4}$ such that $AJ^{-1}(U) \to U$ is a fibration (see \cite{N}, for example), and for a generic fiber, $F_u=AJ^{-1}(u)$, of $AJ$ we have the following short exact sequence: $$\pi_1(F_u) \to \pi_1(Sym_g^\l) \to \pi_1(Pic_g^{4g-4})$$  
Similarly consider the projection $pr:\Q_\l \to \T_g$ given by $(M,q) \mapsto M$. Again, both $\Q_\l$ and $\T_g$ are connected, smooth, complex varieties and $pr$ is a dominant morphism.  Further, notice that $pr^{-1}(M) = AJ^{-1}((M,K_M^2))$. In other words, the fibers of $pr$ may be viewed as a codimension $g$ subset of the fibers of $AJ$. Thus the fibers of $pr$ are also all connected, containing at least one smooth point, and for a generic fiber, $F_M$, of $pr$, we again have the SES: $$\pi_1(F_M) \to \pi_1(\Q_\l) \to \pi_1(\T_g)$$ Since $\T_g$ is simply connected, $\pi_1(\Q_\l)$ is isomorphic to the fundamental group of a generic fiber of $pr$. (One would expect similar arguments apply to any $\l$ such that a generic fiber of $AJ$ is connected - for example any $\l$ of Theorem \ref{thm:c1}). 

If a generic fiber of $pr$ were also a generic fiber of $AJ$ then we would have $\pi_1(\Q_\l) \cong \pi_1(F_u)$ and the kernel of $i_*$ would, in fact, be trivial. However, it is possible that this is not the case, and in general considering the fundamental groups of specific, possibly singular, fibers of $AJ$ seems to be difficult. Again
 for $\l=(1^{4g-4})$, for example, these fibers correspond to the complement of discriminant hypersurfaces in $\P^{3g-4}$ and there is a great deal of literature on the subject, but very little dealing with arbitrarily singular fibers. 

\textbf{Acknowledgments:} I would like to thank P. Seidel, H. Masur, and especially J. Copeland for illuminating conversations and helpful comments.

\end{document}